\documentclass[reqno,a4paper,11pt]{article}
\usepackage{color}
\usepackage[latin1]{inputenc}
\usepackage[english]{babel}
\usepackage{amsmath,amssymb,amsfonts,amsthm,enumerate}
\usepackage{mathrsfs}
\usepackage[T1]{fontenc}
\usepackage{epsfig}
\usepackage{graphicx}
\usepackage[
color,
final] 
{showkeys}

\definecolor{refkeybis}{gray}{.65}
\definecolor{labelkeybis}{gray}{.65}
{\makeatletter
\def\SK@refcolor{\color{refkeybis}}%
\def\SK@labelcolor{\color{labelkeybis}}}

\newtheorem{theorem}{Theorem}[section]
\newtheorem{lemma}[theorem]{Lemma}
\newtheorem{definition}[theorem]{Definition}
\newtheorem{remark}[theorem]{Remark}
\newtheorem{proposition}[theorem]{Proposition}

\newcommand{\Q}{\mathbb{Q}}
\newcommand{\N}{\mathbb{N}}
\newcommand{\R}{\mathbb{R}}

\newcommand{\Leb}[1]{{\mathscr L}^{#1}} 

\renewcommand{\H}{\mathcal{H}}

\renewcommand{\a}{\alpha}

\newcommand{\e}{\varepsilon}
\newcommand{\g}{\gamma}

\renewcommand{\L}{\Lambda}
\newcommand{\n}{\nabla}

\renewcommand{\i}{\infty}
\newcommand{\p}{\partial}

\newcommand{\Cyl}[1]{{\rm Cyl}(#1)}

\renewcommand\div{\operatorname{div}}

\newcommand{\E}{{\mathbb E}}

 \newcommand{\bb}{{\mbox{\boldmath$b$}}}
 \newcommand{\cc}{{\mbox{\boldmath$c$}}}

 \newcommand{\tauV}{{\kern-3pt\tau}}
 
 \newcommand{\sxX}{{\mbox{\scriptsize\boldmath$X$}}}

 \newcommand{\XX}{{\mbox{\boldmath$X$}}}
 \newcommand{\YY}{{\mbox{\boldmath$Y$}}}
 
 \newcommand{\oVVVk}{\overline{\mbox{\boldmath$V$}}\kern-3pt}
 \newcommand{\tVVVk}{\tilde{\mbox{\boldmath$V$}}\kern-3pt}

 \newcommand{\ssigma}{{\mbox{\boldmath$\sigma$}}}

 \newcommand{\eeta}{{\mbox{\boldmath$\eta$}}}

\title{On flows associated to Sobolev vector fields in Wiener spaces: an approach \`a la DiPerna-Lions}
\author{Luigi Ambrosio\
   \thanks{\textsf{l.ambrosio@sns.it}}
   \and
   Alessio Figalli\
   \thanks{\textsf{figalli@unice.fr}}}

\begin{document}

\maketitle

\section{Introduction}

The aim of this paper is the extension to an infinite-dimensional
framework of the theory of flows associated to weakly differentiable
(with respect to the spatial variable $x$) vector fields $\bb(t,x)$.
Starting from the seminal paper \cite{lions}, the finite-dimensional
theory had in recent times many developments, with applications to
fluid dynamics \cite{lions2}, \cite{lions3}, \cite{cufe}, to the
theory of conservation laws \cite{ade}, \cite{adeb}, and it covers
by now Sobolev and even bounded variation \cite{ambrosio}
vectorfields, under suitable bounds on the distributional divergence
of $\bb_t(x):=\bb(t,x)$. Furthermore, in the case of $W^{1,p}_{\rm
loc}$ vector fields with $p>1$, even quantitative error estimates
have been found in \cite{crippade1}; we refer to the Lecture Notes
\cite{cetraro} and \cite{bologna}, and to the bibliographies therein
for the most recent developments on this subject. Our paper fills
the gap, pointed out in \cite{cetraro}, between this family of
results and those available in infinite-dimensional spaces, where
only exponential integrability assumptions on $\nabla\bb_t$ have
been considered so far.

Before passing to the description of our results in Wiener spaces,
we briefly illustrate the heuristic ideas underlying the
above-mentioned finite-dimensional results. The first basic idea is
not to look for pointwise uniqueness statements, but rather to the
family of solutions to the ODE as a whole. This leads to the concept
of flow map $\XX(t,x)$ associated to $\bb$ i.e. a map satisfying
$\XX(0,x)=x$ and $\dot \XX(t,x)=\bb_t(\XX(t,x))$. It is easily seen
that this is not an invariant concept, under modification of $\bb$
in negligible sets. This leads to the concept of $L^r$-regular flow:
we give here the definition adopted in this paper when
$(E,\|\cdot\|)$ is a separable Banach space endowed with a Gaussian
measure $\g$; in the finite-dimensional theory ($E=\R^N$) other
reference measures $\g$ could be considered as well (for instance
the Lebesgue measure \cite{lions}, \cite{ambrosio}).

\begin{definition}[$L^r$-regular $\bb$-flow]\label{dflow1}
Let $\bb:(0,T)\times E\to E$ be a Borel vector field. If
$\XX:[0,T]\times E\to E$ is Borel and $1\leq r\leq\infty$, we say
that $\XX$ is a {\em $L^r$-regular flow} associated to $\bb$ if
the following two conditions hold:
\begin{itemize}
\item[(i)] for $\g$-a.e. $x\in \XX$ the map
$t\mapsto\|\bb_t(\XX(t,x))\|$ belongs to $L^1(0,T)$ and
\begin{equation}\label{bochner}
\XX(t,x)=x+\int_0^t \bb_\tau(\XX(\tau,x))\,d\tau\qquad \forall
t\in [0,T].
\end{equation}
\item[(ii)] for all $t\in [0,T]$ the law of $\XX(t,\cdot)$ under
$\g$ is absolutely continuous with respect to $\g$, with a density
$\rho_t$ in $L^r(\g)$, and $\sup_{t\in [0,T]}\Vert\rho_t\Vert_{L^r(\g)}<\infty$.
\end{itemize}
\end{definition}
In \eqref{bochner}, the integral is understood in Bochner's sense,
namely
$$
\langle e^*,\XX(t,x)-x\rangle=\int_0^t \langle
e^*,\bb_\tau(\XX(\tau,x))\rangle\,d\tau \qquad\forall e^*\in E^*.
$$
It is not hard to show that (see Remark~\ref{invariaflow}),
because of condition (ii), this concept is indeed invariant under
modifications of $\bb$, and so it is appropriate to deal with
vector fields belonging to $L^p$ spaces. On the other hand,
condition (ii) involves all trajectories $\XX(\cdot,x)$ up to
$\g$-neglibigle sets, so the best we can hope for, using this
concept, is existence and uniqueness of $\XX(\cdot,x)$ up to
$\g$-negligible sets.

The second basic idea is the the concept of flow is directly
linked, via the theory of characteristics, to the transport
equation
\begin{equation}\label{transport}
\frac{d}{dt} f(s,x)+\langle\bb_s(x),\nabla_x f(s,x)\rangle=0
\end{equation}
and to the continuity equation
\begin{equation}\label{contieq}
\frac{d}{dt}\mu_t+{\rm div}(\bb_t\mu_t)=0.
\end{equation}
The first link has been exploited in \cite{lions} to transfer
well-posedness results from the transport equation to the ODE,
getting uniqueness of $L^\infty$-regular (with respect to Lebesgue
measure) $\bb$-flows in $\R^N$. This is possible because the flow
maps $(s,x)\mapsto \XX(t,s,x)$ (here we made also explicit the
dependence on the initial time $s$, that we kept equal to $0$ in
Definition~\ref{dflow1}) solve \eqref{transport} for all $t\in
[0,T]$.

Here, in analogy with the approach initiated in \cite{ambrosio}
(see also \cite{figalli} for a stochastic counterpart of it, where
\eqref{contieq} becomes the forward Kolmogorov equation), we
prefer to deal with the continuity equation, which seems to be
more natural in a probabilistic framework. The link between the
ODE and \eqref{contieq} is based on the fact that any positive
finite measure $\eeta$ in $C\bigl([0,T];E\bigr)$ concentrated on
solutions to the ODE is expected to give rise to a weak solution
to \eqref{contieq} (if the divergence operator is properly
understood), with $\mu_t$ given by the marginals of $\eeta$ at
time $t$: indeed, \eqref{contieq} describes the evolution of a
probability density under the action of the ``velocity field''
$\bb$. We shall call these measures $\eeta$ \emph{generalized}
$\bb$-flows. Our goal will be, as in \cite{ambrosio},
\cite{figalli}, to transfer well-posedness informations from the
continuity equation to the ODE, getting existence and uniqueness
results of the $L^r$-regular $\bb$-flows, under suitable
assumptions on $\bb$.

We have to take into account an intrinsic limitation of the theory
of $L^r$-regular $\bb$-flows that is typical of infinite-dimensional
spaces: even if $\bb(t,x)\equiv v$ were constant, the flow map
$\XX(t,x)=x+tv$ would not leave $\g$ quasi-invariant, unless $v$
belongs to a particular subspace of $E$, the so-called
Cameron-Martin space $\H$ of $(E,\g)$, see \eqref{defcm} for its
precise definition. So, from now on we shall assume that $\bb$ takes
its values in $\H$. However, thanks to a suitable change of
variable, we will treat also some non $\H$-valued vector fields, in
the same spirit as in \cite{Peters}, \cite{bogachev}.

We recall that $\H$ can be endowed with a canonical Hilbertian
structure $\langle\cdot,\cdot\rangle_\H$ that makes the inclusion
of $\H$ in $E$ compact; we fix an orthonormal basis $(e_i)$ of
$\H$ and we shall denote by $\bb^i$ the components of $\bb$
relative to this basis (however, all our results are independent
of the choice of $(e_i)$).

With this choice of the range of $\bb$, whenever $\mu_t=u_t\g$ the
equation \eqref{contieq} can be written in the weak sense as
\begin{equation}\label{contieqbis}
\frac{d}{dt}\int_E
u_t\,d\g=\int_E\langle\bb_t,\nabla\phi\rangle_\H u_t\,d\g
\qquad\forall\phi\in\Cyl{E,\g},
\end{equation}
where $\Cyl{E,\g}$ is a suitable space of cylindrical functions
induced by $(e_i)$ (see Definition~\ref{dcyl}). Furthermore, a
Gaussian divergence operator ${\rm div}_\g\cc$ can be defined as
the adjoint in $L^2(\g)$ of the gradient along $\H$:
$$
\int_E\langle\cc,\nabla\phi\rangle_\H\,d\g= -\int_E\phi\,{\rm
div}_\g\cc\,d\g\qquad\forall\phi\in\Cyl{E,\g}.
$$
Another typical feature of our Gaussian framework is that
$L^\infty$-bounds on ${\rm div}_\g$ do not seem natural, unlike
those on the Euclidean divergence in $\R^N$ when the reference
measure is the Lebesgue measure: indeed, even if
$\bb(t,x)=\cc(x)$, with $\cc:\R^N\to\R^N$ smooth and with bounded
derivatives, we have ${\rm div}_\g\cc={\rm
div}\cc-\langle\cc,x\rangle$ which is unbounded, but exponentially
integrable with respect to $\g$.

We can now state the main result of this paper:

\begin{theorem}[Existence and uniqueness of $L^r$-regular $\bb$-flows]\label{main1}
Let $p,\,q>1$ and let $\bb:(0,T)\times E\to\H$ be satisfying:
\begin{itemize}
\item[(i)] $\|\bb_t\|_\H\in L^1\bigl((0,T);L^p(\g)\bigr)$;
\item[(ii)] for a.e. $t\in (0,T)$ we have $\bb_t\in
LD^q_\H(\g;\H)$ with
\begin{equation}\label{ldsym}
\int_0^T\biggl(\int_E \|(\nabla\bb_t)^{\rm
sym}(x)\|^q_{HS}\,d\g(x)\biggr)^{1/q}\, dt<\infty,
\end{equation}
and ${\rm div}_\g\bb_t\in L^1\bigl((0,T);L^q(\g)\bigr)$;
\item[(iii)] $\exp (c[{\rm div}_\g \bb_t]^-)\in
L^\infty\bigl((0,T);L^1(\g)\bigr)$ for some $c>0$.
\end{itemize}
If $r:=\max\{p',q'\}$ and $c\geq rT$, then the $L^r$-regular flow
exists and is unique in the following sense: any two $L^r$-regular
flows $\XX$ and $\tilde \XX$ satisfy
$$
\XX(\cdot,x)=\tilde{\XX}(\cdot,x)\quad\text{in $[0,T]$, for
$\g$-a.e. $x\in E$.}
$$
Furthermore, $\XX$ is $L^s$-regular for all $s\in [1,\frac{c}{T}]$
and the density $u_t$ of the law of $\XX(t,\cdot)$ under $\g$
satisfies
$$
\int (u_t)^s\,d\g\leq \biggl\|\int_E\exp\bigl(Ts[{\rm
div}_\g\bb_t]^-\bigr)\,d\g\biggr\Vert_{L^\infty(0,T)}
\qquad\text{for all $s\in [1,\frac{c}{T}]$.}
$$
In particular, if $\exp (c[{\rm div}_\g \bb_t]^-)\in
L^\infty\bigl((0,T);L^1(\g)\bigr)$ for all $c>0$, then the
$L^r$-regular flow exists globally in time, and is $L^s$-regular for
all $s \in [1,\i).$
\end{theorem}

The symmetric matrix $(\nabla\bb_t)^{\rm sym}$, whose
Hilbert-Schmidt norm appears in \eqref{ldsym}, corresponds to the
symmetric part of the derivative of $\bb_t$, defined in a weak
sense by \eqref{explsymm}: notice that, in analogy with the finite
dimensional result \cite{capuzzo}, \emph{no} condition is imposed
on the antisymmetric part of the derivative, which need not be
given by a function; this leads to a particular function space
$LD^q(\g;\H)$ (well studied in linear elasticity in finite
dimensions, see \cite{temam}) which is for instance larger than
the Sobolev space $W^{1,q}_{\H}(\g;\H)$, see
Definitions~\ref{defSobolev} and \ref{defLD}. Also, we will prove
that uniqueness of $\XX$ holds even within the larger class of
generalized $\bb$-flows.

Let us explain first the main differences between our strategy and
the techniques used in \cite{cruzeiro1}, \cite{cruzeiro2},
\cite{cruzeiro3}, \cite{Peters}, \cite{bogachev} for autonomous
(i.e. time independent) vector fields in infinite-dimensional
spaces. The standard approach for the existence of a flow consists
in approximating the vector field $\bb$ with finite-dimensional
vector fields $\bb_N$, constructing a finite-dimensional flow
$\XX_N$, and then passing to the limit as $N\to\infty$. This part
of the proof requires quite strong a-priori estimates on the flows
to have enough compactness to pass to the limit. To get these
a-priori estimates, the assumptions on the vector field, instead
of the hypotheses (i)-(iii) in Theorem~\ref{main1}, are:
$$
\|\bb\|_\H \in \bigcap_{p\in [1,\infty)} L^p(\g),
$$
$$
\exp (c\|\n \bb\|_{{\cal L}(\H,\H)})\in L^1(\g) \qquad \text{for
all }c>0,
$$
$$
\exp (c|{\rm div}_\g \bb|)\in L^1(\g) \qquad \text{for some }c>0,
$$
where $\|\n \bb\|_{{\cal L}(\H,\H)}$ denotes the operator norm of
$\n \bb$ from $\H$ to $\H$. So, apart from the minor fact that we
allow a measurable time dependence of $\bb$, the main difference
between these results and ours is that we replace exponential
integrability of $\bb$ and the \emph{operator} norm of $\nabla\bb$ by $p$-integrability of
$\bb$ and $q$-integrability of the (stronger) \emph{Hilbert-Schmidt} norm of
$\nabla\bb_t$ (or, as we said, of its symmetric part).

Let us remark for instance that, just for the existence part of a
generalized $\bb$-flow, the hypothesis on $\div_\g\bb$ could be
relaxed to a one sided bound, as we did. Indeed, this assumption
allows to prove uniform estimates on the density of the
approximating flows, see for instance Theorem~\ref{texflow}. On
the other hand, the proof of the uniqueness of the flow strongly
relies on the fact that one can use the approximating flows
$\XX_N$ also for negative times.

Our strategy is quite different from the above one: the existence
and uniqueness of a regular flow will be proved at once in the
following way. First of all, the existence of a generalized
$\bb$-flow $\eeta$, even without the regularity assumption
\eqref{ldsym}, can be obtained thanks to a tightness argument for
measures in $C\bigl([0,T];E\bigr)$ and proving uniform estimates on
the density of the finite-dimensional approximating flows. Then we
prove uniqueness in the class of generalized $\bb$-flows. This
implies as a byproduct that $\eeta$ is induced by a
``deterministic'' $\XX$, thus providing the desired existence and
uniqueness result. Moreover the flexibility of this approach allows
us to prove the stability of the $L^r$-regular flow under smooth
approximations of the vector field, and thanks to the uniqueness we
can also easily deduce the semigroup property.

The main part of the paper is therefore devoted to the proof of
uniqueness. As we already said, this depends on the well-posedness
of the continuity equation \eqref{contieqbis}. Specifically, we
will show uniqueness of solutions $u_t$ in the class
$L^\infty\bigl((0,T);L^r(\g)\bigr)$. The key point, as in the
finite-dimensional theory, is to pass from \eqref{contieqbis} to
\begin{equation}\label{contieqter}
\frac{d}{dt}\int_E
\beta(u_t)\,d\g=\int_E\langle\bb_t,\nabla\phi\rangle_\H
\beta(u_t)\,d\g+ \int_E [\beta(u_t)-u_t\beta'(u_t)]{\rm
div}_\g\bb_t\,d\g \qquad\forall\phi\in\Cyl{E,\g},
\end{equation}
for all $\beta\in C^1(\R)$ with $\beta'(z)$ and
$z\beta'(z)-\beta(z)$ bounded, and then to choose as function
$\beta$ suitable $C^1$ approximations of the positive or of the
negative part, to show that the equation preserves the sign of the
initial condition. The passage from \eqref{contieqbis} to
\eqref{contieqter} can be \emph{formally} justified using the rule
$$
{\rm div}_\gamma (v\cc)=v{\rm div}_\gamma\cc+\langle\nabla
v,\cc\rangle_\H
$$
and the chain rule $\n \beta(u)=\beta'(u)\n u$, but it is not
always possible. It is precisely at this place that the regularity
assumptions on $\bb_t$ enter. The finite-dimensional strategy
involves a regularization argument (in the space variable only)
and a careful analysis of the ``commutators'' (with $v=u_t$,
$\cc=\bb_t$)
$$
r^\e(\cc,v):=e^\e\langle\cc,\nabla T_\e v\rangle_\H-T_\e({\rm
div}_\g(v\cc)),
$$
where $\e$ is the regularization parameter and $T_\e$ is the
regularizing operator. Already in the finite-dimensional theory (see
\cite{lions}, \cite{ambrosio}) a careful estimate of $r^\e$ is
needed, taking into account some cancellation effects. These effects
become even more important in this framework, where we use as a
regularizing operator the Ornstein-Uhlenbeck operator \eqref{Mehler}
(in particular the semigroup property and the fact that $T_t$ is
self-adjoint from $L^p(\g)$ to $L^{p'}(\g)$ will play an important
role). The core of our proof is indeed Section~\ref{sect:commu},
where we obtain commutator estimates in $\R^N$ independent of $N$,
and therefore suitable for an extension, via the canonical
cylindrical approximation, to $E$.

The paper is structured as follows: first we recall the main
notation needed in the paper. In Section~\ref{sect:wellposed} we
prove the well-posedness of the continuity equation, while in
Section~\ref{sect:ODE} we prove existence, uniqueness and
stability of regular flows. The results of both sections rely on
some finite dimensional a-priori estimates that we postpone to
Section~\ref{sect:finitedim}. Finally, to apply our results also
in more general situations: in Section~\ref{sect:nonCM} we see how
our results can be extended to the case non $\H$-valued vector
fields, in the same spirit as in \cite{Peters}, \cite{bogachev}.

\section{Main notation and preliminary results}

\smallskip
{\bf Measure-theoretic notation.} All measures considered in this
paper are positive, finite and defined on the Borel
$\sigma$-algebra. Given $f:E\to F$ Borel and a measure $\mu$ in
$E$, we denote by $f_\#\mu$ the push-forward measure in $F$, i.e.
the law of $f$ under $\mu$. We denote by $\chi_A$ the
characteristic function of a set $A$, equal to 1 on $A$, and equal
to 0 on its complement.

We consider a separable Banach space $(E,\|\cdot\|)$ endowed with
a Gaussian measure $\g$, i.e. $(e^*)_\#\g$ is a Gaussian measure
in $\R$ for all $e^*\in E^*$. We shall assume that $\g$ is
centered and non-degenerate, i.e. that $\int_E x\,d\g(x)=0$ and
$\g$ is not supported in a proper subspace of $E$. We recall (see
\cite{Ledoux}) that, by Fernique's theorem, $\int_E \exp(
c\|x\|^2)\,d\g(x)<\infty$, whenever $2c<\sup_{\|e^*\|\leq 1}
\|\langle e^*,x\rangle\|_{L^2(\g)}$.

\smallskip
\noindent {\bf Cameron-Martin space.} We shall denote by $\H\subset
E$ the Cameron Martin space associated to $(E,\g)$. It can be
defined \cite{Bogachev,Ledoux} as
\begin{equation}\label{defcm}
\H:=\left\{\int_E \phi(x)x\,d\g(x)\,:\,\phi\in L^2(\g)\right\}.
\end{equation}
The non-degeneracy assumption assumption on $\g$ easily implies
that $\H$ is a dense subset of $E$. If we denote by $i:L^2(\g)\to
\H\subset E$ the map $\phi\mapsto\int_E\phi(x)x\,d\g(x)$, and by
$K$ the kernel of $i$, we can define the Cameron-Martin norm
$$
\|i(\phi)\|_\H:=\min_{\psi\in K}\|\phi-\psi\|_{L^2(\g)},
$$
whose induced scalar product $\langle\cdot,\cdot\rangle_\H$
satisfies
\begin{equation}\label{repscaH}
\langle i(\phi),i(\psi)\rangle_\H=\int_E\phi\psi\,d\g\qquad
\forall\phi\in L^2(\g),\,\,\forall\psi\in K^\perp.
\end{equation}
Notice also that $i(\langle e^*,x\rangle)\in K^\perp$ for all
$e^*\in E^*$, because
$$
\int_E \langle e^*,x\rangle\psi(x)\,d\g(x)=\langle e^*,\int_E
x\psi(x)\,d\g(x)\rangle=0\qquad\forall\psi\in K.
$$
Since $i$ is not injective in general, it is often more convenient
to work with the map $j:E^*\to\H$, dual of the inclusion map of $\H$
in $E$ (i.e. $j(e^*)$ is defined by $\langle
j(e^*),h\rangle_\H=\langle e^*,h\rangle$ for all $h\in\H$). The set
$j(E^*)$ is obviously dense in $\H$ (for the norm $\|\cdot\|_\H$),
and $j$ is injective thanks to the density of $\H$ in $E$;
furthermore, choosing $\phi(x)=\langle e^*,x\rangle$ in
\eqref{repscaH}, we see that $i(\langle e^*,x\rangle)=j(e^*)$. As a
consequence the vector space $\{\langle e^*,x\rangle\,:\,e^*\in
E^*\}$ is dense in $K^\perp$. Since $\|i(\langle
e^*,x\rangle)\|\leq\bigl(\int_E\|x\|^2\,d\g\bigr)^{1/2}\Vert\langle
e^*,x\rangle\Vert_{L^2(\g)}=\|i(\langle e^*,x\rangle)\|_\H$, the
inclusion of $\H$ in $E$ is continuous, and it is not hard to show
that it is also compact (see \cite[Corollary~3.2.4]{Bogachev}).

This setup becomes much simpler when $(E,\|\cdot\|)$ is an Hilbert
space:

\begin{remark} [The Hilbert case]\label{rhilbert} {\rm
Assume that $(E,\|\cdot\|)$ is an Hilbert space. Then, after
choosing an orthonormal basis in which the covariance operator
$(x,y)\mapsto\int_E\langle x,z\rangle\langle y,z\rangle\,d\g(z)$
is diagonal, we can identify $E$ with $\ell^2$, endowed with the
canonical basis $\epsilon_i$, and the coordinates $x_i$ of
$x\in\ell^2$ relative to $\epsilon_i$ are independent, Gaussian
and with variance $\lambda_i^2$ (with $\lambda_i>0$ by the
non-degeneracy assumption). Then, the integrability of $\|x\|^2$
implies that $\sum_i\lambda_i^2$ is convergent, $e_i^*=\epsilon_i$
(here we are using the Riesz isomorphism to identify $\ell^2$ with
its dual), $e_i=\lambda_i\epsilon_i$ and the Cameron-Martin space
is
$$
\H:=\left\{x\in\ell^2\,:\,\sum_{i=1}^\infty
\frac{(x^i)^2}{\lambda_i^2}<\infty\right\}.
$$
The map $j:\ell^2\to\H$ is given by $(x_i)\mapsto (\lambda_ix_i)$.
}\end{remark}

Let us remark that, although we constructed $\H$ starting from $E$,
it is indeed $\H$ which plays a central role in our results;
according to the Gross viewpoint, this space might have been taken
as the starting point, see \cite[\S3.9]{Bogachev} and
Section~\ref{finitedim} for a discussion of this fact.

\smallskip
\noindent {\bf Finite-dimensional projections.} The
above-mentioned properties of $j$ allow the choice of
$(e_n^*)\subset E^*$ such that $(j(e_n^*))$ is a complete
orthonormal system in $\H$. Then, setting $e_n:=j(e_n^*)$, we can
define the continuous linear projections $\pi_N:E\to\H$ by
\begin{equation}\label{defpiN}
\pi_N(x):=\sum_{k=1}^N \langle e_k^*,x\rangle
e_k\biggl(=\sum_{k=1}^N\langle e_k,x\rangle_\H e_k\quad\text{for
$x\in\H$}\biggr).
\end{equation}
The term ``projection'' is justified by the fact that, by the
second equality in \eqref{defpiN}, $\pi_N\vert_\H$ is indeed the
orthogonal projection on
\begin{equation}\label{defHN}
\H_N:={\rm span\,}\bigl(e_1,\ldots,e_N\bigr).
\end{equation}
From now such a basis $(e_i)$ of $\H$ will be fixed, and we shall
denote by $v^i$ the components of $v\in\H$ relative to this basis.
Also, for a given Borel function $u:E\to\R$, we shall denote by
$\E_Nu$ the conditional expectation of $u$ relative to the
$\sigma$-algebra generated by $\langle e_1^*,x\rangle,\ldots,
\langle e_N^*,x\rangle$. The following result follows by martingale
convergence theorems, because the $\sigma$--algebra generated by
$\langle e_i^*,x\rangle$ is the Borel $\sigma$-algebra (see also
\cite[Corollary~3.5.2]{Bogachev}):

\begin{lemma}\label{lprod}
For all $p\in [1,\infty)$ and $u\in L^p(\g)$ we have $\E_Nu\to u$
$\g$-a.e. and in $L^p(\g)$.
\end{lemma}

According to these projections, we can define the space
$\Cyl{E,\g}$ of smooth cylindrical functions (notice that this
definition depends on the choice of the basis $(e_n)$).

\begin{definition}[Smooth cylindrical functions]\label{dcyl}
Let $C^\infty_b(\R^N)$ be the space of smooth functions in $\R^N$,
bounded together with all their derivatives. We say that
$\phi:E\to\R$ is {\em cylindrical} if
\begin{equation}\label{repcyl}
\phi(x)=\psi\bigl(\langle e_1^*,x\rangle,\ldots,\langle
e_N^*,x\rangle\bigr)
\end{equation}
for some integer $N$ and some $\psi\in C^\infty_b(\R^N)$.
\end{definition}

If $v\in E$ and $\phi:E\to\R$ we shall denote by $\partial_v\phi$
the partial derivative of $\phi$ along $v$, wherever this exists.
Obviously, cylindrical functions are differentiable infinitely
many times in all directions: if $\phi$ is as in \eqref{repcyl},
the first order derivative is given by
\begin{equation}\label{svitla}
\partial_v\phi(x)=\sum_{i=1}^N\frac{\partial\psi}{\partial z_i}
\bigl(\langle e_1^*,x\rangle,\ldots,\langle e_N^*,x\rangle\bigr)
\langle e_i^*,v\rangle.
\end{equation}
If $v\in\H$ the above formula becomes
$$
\partial_v\phi(x)=\sum_{i=1}^N\frac{\partial\psi}{\partial z_i}
\bigl(\langle e_1^*,x\rangle,\ldots,\langle e_N^*,x\rangle\bigr)
\langle e_i,v\rangle_{\H},
$$
and this allows to define the gradient of $\phi$ as an element of
$\H$:
$$
\n \phi(x):=\sum_{i=1}^N\frac{\partial\psi}{\partial z_i}
\bigl(\langle e_1^*,x\rangle,\ldots,\langle e_N^*,x\rangle\bigr)
e_i \in \H.
$$

\smallskip
\noindent {\bf Gaussian divergence and differentiability along
$\H$.} Let $\bb:E\to\H$ be a vector field with $\|\bb\|_\H\in
L^1(\g)$; we say that a function ${\rm div}_\g \bb\in L^1(\g)$ is
the \emph{Gaussian divergence} of $b$ (see for instance
\cite[\S5.8]{Bogachev}) if
\begin{equation}\label{defgdiv}
\int_E\langle\nabla\phi,\bb\rangle_\H\,d\g=-\int_E \phi\,{\rm
div}_\g \bb\,d\g \qquad\forall\phi\in\Cyl{E,\g}.
\end{equation}
In the finite-dimensional space $E=\R^N$ endowed with the standard
Gaussian we have, by an integration by parts,
\begin{equation}\label{expdiv}
{\rm div}_\g\bb={\rm div\,}\bb-\langle \bb,x\rangle.
\end{equation}

We recall the integration by parts formula
\begin{equation}\label{basic_integration}
\int_E \partial_{j(e^*)}\phi\,d\g=\int_E \phi \langle
e^*,x\rangle\,d\g \qquad\forall\phi\in\Cyl{E,\g},\,\, \forall
e^*\in E^*.
\end{equation}
This motivates the following definitions: if both $u(x)$ and
$u(x)\langle e^*,x\rangle$ belong to $L^1(\g)$, we call weak
derivative of $u$ along $j(e^*)$ the linear functional on
$\Cyl{E,\g}$
\begin{equation}\label{veryweak}
\phi\mapsto -\int_E u\partial_{j(e^*)}\phi\,d\g+\int_E
u\phi\langle e^*,x\rangle\,d\g.
\end{equation}
As in the classical finite-dimensional theory, we can define Sobolev
spaces by requiring that these functionals are representable by
$L^q(\g)$ functions, see Chapter~5 of \cite{Bogachev} for a more
complete discussion of this topic.

\begin{definition}[Sobolev space $W^{1,q}_\H(\g)$]
\label{defSobolev} If $1\leq q\leq\infty$, we say that $u\in
L^1(\g)$ belongs to $W^{1,q}_{\H}(E,\g)$ if $u(x)\langle
e^*,x\rangle\in L^1(\g)$ for all $e^*\in E^*$ and there exists
$g\in L^q(\g;\H)$ satisfying
\begin{equation}\label{sciopero}
\int_E u\partial_{j(e^*)}\phi\,d\g+\int_E\phi\langle
g,j(e^*)\rangle_\H\,d\g=\int_E u\phi\langle e^*,x\rangle\,d\g
\qquad\quad\forall e^*\in E^*,\,\,\forall\phi\in\Cyl{E,\g}.
\end{equation}
\end{definition}

The condition $u(x)\langle e^*,x\rangle\in L^1(\g)$ is
automatically satisfied whenever $u\in L^p(\g)$ for some $p>1$,
thanks to the fact that the law of $\langle e^*,x\rangle$ under
$\g$ is Gaussian, so that $\langle e^*,x\rangle\in L^r(\g)$ for
all $r<\infty$.

We shall denote, as usual, the (unique) weak derivative $g$ by
$\nabla u$ and its components $\langle g,e_i\rangle_\H$ by
$\partial_i u$, so that \eqref{sciopero} becomes
\begin{equation}\label{sciopero1}
\int_E u\partial_i\phi\,d\g+\int_E\phi\partial_i u\,d\g=\int_E
u\phi\langle e^*_i,x\rangle\,d\g \qquad\quad\forall i\geq
1,\,\,\forall\phi\in\Cyl{E,\g}.
\end{equation}

We recall that a continuous linear operator $L:\H\to\H$ is said to
be Hilbert-Schmidt if $\|L\|_{HS}$, defined as the square root of
the trace of $L^tL$, is finite. Accordingly, if $L_{ij}=\langle
L(e_i),e_j\rangle_\H$ is the symmetric matrix representing
$L:\H\to\H$ in the basis $(e_i)$, we have that $L$ is of
Hilbert-Schmidt class if and only if $\sum_{ij}L_{ij}^2$ is
convergent, and
\begin{equation}\label{defHS}
\|L\|_{HS}=\sqrt{\sum_{ij}L^2_{ij}}.
\end{equation}

The following proposition shows that bounded continuous operators
from $E$ to $\H$ are of Hilbert-Schmidt class, when restricted to
$\H$. In particular our results apply under $p$-integrability
assumptions on $\nabla\bb_t$ when the operator norm between $E$ and
$\H$ is used.

\begin{proposition}
\label{propCauchy} Let $L:E\to\H$ be a linear continuous operator.
Then the restriction of $L$ to $\H$ is of Hilbert-Schmidt class and
$\|L\|_{HS}\leq C\|L\|_{{\cal L}(E,\H)}$, with $C$ depending only on
$E$ and $\g$.
\end{proposition}
\begin{proof}
By \cite[Theorem 3.5.10]{Bogachev} we can find a complete
orthonormal system $(f_n)$ of $\H$ such that
$\sum_n\|f_n\|^2=:C<+\infty$. Denoting by $\|L\|$ the operator norm
of $L$ from $E$ to $\H$, we have then
$$
\|L\|_{HS}^2=\sum_{i,j}(\langle
L(f_i),f_j\rangle_\H)^2=\sum_i\||L(f_i)\|_\H^2\leq\|L\|^2 \sum_i
\|f_i\|^2=C\|L\|^2.
$$
\end{proof}

From now on, we shall denote by $L^p(\g;\H)$ the space of Borel
maps $\cc:E\to\H$ such that $\|\cc\|_\H\in L^p(\g)$. Given the
basis $(e_i)$ of $\H$, we shall denote by $\cc^i$ the components
of $\cc$ relative to this basis.

\begin{definition}[The space $LD(\g;\H)$]\label{defLD}
If $1\leq q\leq\infty$, we say that $\cc\in L^1(\g;\H)$ belongs to
$LD^q(\g;\H)$ if:
\begin{itemize}
\item[(a)] for all $h=j(e^*)\in\H$, the function
$\langle\cc,h\rangle_\H$ has a weak derivative in $L^q(\g)$ along
$h$, that we shall denote by $\partial_h\langle\cc,h\rangle_\H$,
namely
\begin{equation}\label{sciopero2}
\int_E \langle\cc,h\rangle_\H\partial_h\phi\,d\g+\int_E\phi
\partial_h\langle\cc,h\rangle_\H\,d\g=\int_E \langle\cc,h\rangle_\H\phi\langle
e^*,x\rangle\,d\g \qquad\quad\forall\phi\in\Cyl{E,\g};
\end{equation}
\item[(b)] the symmetric matrices
\begin{equation}\label{defsymmo}
(\nabla\cc)_{ij}^{\rm sym}(x):=\frac{1}{4}
\bigl[\partial_{(e_i+e_j)}(\cc^i+\cc^j)(x)-\partial_{(e_i-e_j)}(\cc^i-\cc^j)(x)\bigr]
\end{equation}
satisfy
$$
\int_E\|(\nabla\cc)^{\rm sym}\|^q_{HS}\,d\g<\infty.
$$
\end{itemize}
\end{definition}

If all components $\cc^i$ of $\cc$ belongs to $W^{1,q}_\H(\g)$
then the function $(\nabla\cc)^{\rm sym}_{ij}$ in \eqref{defsymmo}
really corresponds to the symmetric part of $(\nabla
c)_{ij}=\partial_j\cc^i$, and this explains our choice of
notation. However, according to our definition of $LD^q(\g;\H)$,
the vector fields $\cc$ in this space need not have components
$\cc^i$ in $W^{1,q}_\H(\g)$. Moreover, from \eqref{defsymmo} we
obtain that $(\partial_i\cc^j+\partial_j\cc^i)/2$ are
representable by the $L^q(\g)$ functions $(\nabla\cc)^{\rm
sym}_{ij}$, namely
\begin{equation}\label{explsymm}
\int_E
\frac{1}{2}(\cc^i\partial_j\phi+\cc^j\partial_i\phi)\,d\g+\int_E
\phi (\nabla\cc)_{ij}^{\rm
sym}\,d\g=\int_E\frac{1}{2}(\cc^i\langle
e_j^*,x\rangle+\cc^j\langle e_i^*,x\rangle)\phi\,d\g \qquad\forall
\phi\in\Cyl{E,\g}.
\end{equation}

\begin{remark}[Density of cylindrical functions]\label{denscyl}{\rm
We recall that $\Cyl{E,\g}$ is dense in all spaces $W^{1,p}_\H(\g)$,
$1\leq p<\infty$. More precisely, if $1\leq p,q<\infty$, any
function $u\in W^{1,p}_\H(\g)\cap L^q(\g)$ can be approximated in
$L^q(\g)$ by cylindrical functions $u_n$ with $\nabla u_n\to\nabla
u$ strongly in $L^p(\g;\H)$. In the case $p=\infty$, convergence of
the gradients occurs in the weak$^*$ topology of $L^\infty(\g;\H)$.
These density results can be proved first in the finite-dimensional
case and then, thanks to Lemma~\ref{lprod}, in the general
case.}\end{remark}

\begin{remark}\label{rprodu} {\rm
In the sequel we shall use the simple rule
$$
{\rm div}_\g(\bb u)=u{\rm div}_\g\bb+\langle\bb,\nabla
u\rangle_\H,
$$
valid whenever ${\rm div}_\g\bb\in L^p(\g)$, $u\in L^{p'}(\g)$,
$\bb\in L^q(\g;\H)$ and $u\in W^{1,q'}_\H(\g)$. The proof is a
direct consequence of Remark~\ref{denscyl}. }\end{remark}

\begin{remark} [Invariance of $\div_\g$, $W^{1,q}_\H(\g)$,
$LD^q(\g)$]\label{rinvaria} {\rm The definitions of Gaussian
divergence, Sobolev space and $LD$ space, as given, involve the
space $\Cyl{E,\g}$, which depends on the choice of the complete
orthonormal basis $(e_i)$. However, an equivalent formulation could
be given using the space $C^1_b(E,\g)$ of functions that are Frechet
differentiable along all directions in $\H$, with a bounded
continuous gradient: indeed, cylindrical functions belong to
$C^1_b(E,\g)$, and since $C^1_b(E,\g)$ is contained in
$W^{1,\infty}_\H(\g)$, thanks to Remark~\ref{denscyl} the functions
in this space can be well approximated (in all spaces $L^p(\g)$ with
$p<\infty$, and with weak$^*$ convergence in $L^\infty(\g)$ of
gradients) by cylindrical functions. A similar remark applies to the
continuity equation, discussed in the next section.}\end{remark}

\section{Well posedness of the continuity equation}
\label{sect:wellposed}

Let $I\subset\R$ be an open interval. In this section we shall
consider the \emph{continuity equation} in $I\times E$, possibly
with a source term $f$, i.e.
\begin{equation}\label{conti}
\frac{d}{dt}(u_t\g) +{\rm div}_\g(\bb_t u_t\g)=f\g.
\end{equation}
This equation has to be understood in the weak sense, namely we
require that $t\mapsto\int_E u_t\phi\,d\g$ is absolutely
continuous in $I$ and
\begin{equation}\label{conti2}
\frac{d}{dt}\int_E u_t\phi\,d\g=\int_E\langle
\bb_t,\nabla\phi\rangle_\H u_t\,d\g+ \int_E f\phi\,d\g
\qquad\text{a.e. in $I$, for all $\phi\in\Cyl{E,\g}$.}
\end{equation}
The minimal requirement necessary to give a meaning to
\eqref{conti2} is that $u$, $f$ and $|u|\|\bb\|_\H$ belong to
$L^1\bigl(I;L^1(\g)\bigr)$, and we shall always make assumptions
on $u$, $f$ and $\bb$ to ensure that these properties are
satisfied.

Sometimes, to simplify our notation, with a slight abuse we drop
$\g$ and write \eqref{conti} just as
$$
\frac{d}{dt}u_t+{\rm div}_\g(\bb_t u_t)=f.
$$
However, we always have in mind the weak formulation
\eqref{conti2}, and we shall always assume that $f\in
L^1\bigl(I;L^1(\g)\bigr)$.

Since we are, in particular, requiring all maps $t\mapsto \int_E
u_t\phi\,d\g$ to be uniformly continuous in $I$, the map $t\mapsto
u_t$ is weakly continuous in $I$, with respect to the duality of
$L^1(\g)$ with $\Cyl{E,\g}$. Therefore, if $I=(0,T)$, it makes
sense to say that a solution $u_t$ of the continuity equation
starts from $\bar u\in L^1(\g)$ at $t=0$:
\begin{equation}\label{initialcond}
\lim_{t\downarrow 0}\int_E u_t\phi\,d\g=\int_E\bar
u\phi\,d\g\qquad\forall u\in\Cyl{E,\g}.
\end{equation}

\begin{theorem}[Well-posedness of the continuity equation]\label{main2}
(Existence) Let $\bb:(0,T)\times E\to\H$ be satisfying
\begin{equation}\label{intb}
\|\bb_t\|_\H\in L^1\bigl((0,T);L^p(\g)\bigr)\qquad\text{for some
$p>1$}
\end{equation}
and
\begin{equation}\label{peeters}
\exp(c[{\rm div}_\g \bb_t]^-)\in
L^\infty\bigl((0,T);L^1(\g)\bigr)\qquad\text{for some $c>T p'$.}
\end{equation}
Then, for any nonnegative $\bar u\in L^\i(\g)$, the continuity
equation has a nonnegative solution $u_t$ with $u_0=\bar u$
satisfying (as a byproduct of its construction)
\begin{equation}\label{expbound}
\int (u_t)^r\,d\g\leq \|\bar u\|_{L^\i(\g)}^r \biggl\|\int_E
\exp\bigl(Tr [{\rm
div}_\g\bb_t]^-\bigr)\,d\g\biggr\Vert_{L^\infty(0,T)} \qquad
\text{for all } r \in [1,\frac{c}{T}], \,t \in [0,T].
\end{equation}
\noindent (Uniqueness) Let $\bb:(0,T)\times E\to\H$ be satisfying
\eqref{intb}, $\bb_t\in LD^q(\g;\H)$ for a.e. $t\in (0,T)$ with
\begin{equation}\label{HS}
\int_0^T\biggl(\int_E \|(\nabla\bb_t)^{\rm
sym}\|^q_{HS}\,d\g\biggr)^{1/q}\,dt<\infty
\end{equation}
for some $q>1$, and
\begin{equation}\label{divbe}
{\rm div}_\g\bb_t\in L^1\bigl((0,T);L^q(\g)\bigr).
\end{equation}
Then, setting $r=\max\{p',q'\}$, if $c\geq Tr$ the continuity
equation \eqref{conti} in $(0,T)\times E$ has at most one solution
in the function space $L^\infty\bigl((0,T);L^r(\g)\bigr)$.
\end{theorem}

\begin{definition}[Renormalized solutions] We say that a solution
$u_t$ of \eqref{conti} in $I\times E$ is {\em renormalized} if
\begin{equation}\label{conti1}
\frac{d}{dt}\beta(u_t)+{\rm div}_\g(\bb_t\beta(u_t))=
[\beta(u_t)-u_t\beta'(u_t)]{\rm div}_\g \bb_t+ f\beta'(u_t)
\end{equation}
in the sense of distributions in $I\times E$, for all $\beta\in
C^1(\R)$ with $\beta'(z)$ and $z\beta'(z)-\beta(z)$ bounded.
\end{definition}

In the sequel we shall often use the Ornstein-Uhlenbeck operator
$T_t$, defined for $u\in L^1(\g)$ by Mehler's formula
\begin{equation}\label{Mehler}
T_t u(x):=\int_E u(e^{-t}x+\sqrt{1-e^{-2t}}y)\,d\g(y).
\end{equation}

In the next proposition we summarize the main properties of the OU
operator used in this paper, see Theorem~1.4.1, Theorem~2.9.1 and
Proposition~5.4.8 of \cite{Bogachev}.

\begin{proposition}[Properties of the OU semigroup]\label{pOU}
Let $T_t$ be as in \eqref{Mehler}.
\begin{itemize}
\item[(i)] $\|T_t u\|_{L^p(\g)}\leq\|u\|_{L^p(\g)}$ for all $u\in
L^p(\g)$, $p\in [1,\infty]$, $t\geq 0$, and equality holds if $u$ is
nonnegative and $p=1$.
\item[(ii)] $T_t$ is self-adjoint in $L^2(\g)$ for all $t\geq 0$.
More generally, if $1\leq p\leq\infty$, we have
\begin{equation}\label{selfdual}
\int_E v T_t u\,d\g=\int_E u T_t v\,d\g\qquad \forall u\in
L^p(\g),\,\,\forall v\in L^{p'}(\g).
\end{equation}
\item[(iii)] For all $p\in (1,\infty)$, $t>0$ and $u\in L^p(\g)$
we have $T_t u\in W^{1,p}_\H(\g)$ and
\begin{equation}
\|\nabla T_t u\|_{L^p(\g;\H)}\leq C(p,t)\|u\|_{L^p(\g)}.
\end{equation}
\item[(iv)] For all $p\in [1,\infty]$ and $u\in W^{1,p}_\H(\g)$ we
have $\nabla T_t u= e^{-t}T_t\nabla u$.
\item[(v)] $T_t$ maps
$\Cyl{E,\g}$ in $\Cyl{E,\g}$ and $T_t u\to u$ in $L^p(\g)$ as
$t\downarrow 0$ for all $u\in L^p(\g)$, $1\leq p<\infty$.
\end{itemize}
\end{proposition}

In the same spirit of \eqref{veryweak}, we can now extend the
action of the semigroup from $L^1(\g)$ to elements $\ell$ in the
algebraic dual of $\Cyl{E,\g}$ as follows:
$$
\langle
T_t\ell,\phi\rangle:=\langle\ell,T_t\phi\rangle\qquad\phi\in\Cyl{E,\g}.
$$
This is an extension, because if $\ell$ is induced by some
function $u\in L^1(\g)$, i.e. $\langle\ell,\phi\rangle=\int_E\phi
u\,d\g$ for all $\phi\in\Cyl{E,\g}$, then because of
\eqref{selfdual} $T_t\ell$ is induced by $T_t u$, i.e. $\langle
T_t\ell,\phi\rangle=\int_E \phi T_t u\,d\g$ for all
$\phi\in\Cyl{E,\g}$. In general we shall say that $T_t\ell$ is a
function whenever there exists (a unique) $v\in L^1(\g)$ such that
$\langle T_t\ell,\phi\rangle=\int_Ev\phi\,d\g$ for all $\phi\in
\Cyl{E,\g}$.

In the next lemma we will use this concept when $\ell$ is the
Gaussian divergence of a vector field $\cc$: indeed, $\ell$ can be
thought via the formula
$-\int_E\langle\cc,\nabla\phi\rangle_\H\,d\g$ as an element of the
dual of $\Cyl{E,\g}$. Our first proposition provides a sufficient
condition ensuring that $T_t({\rm div}_\g\cc)$ is a function.

\begin{lemma}\label{hyperco}
Assume that $r\in (1,\infty)$ and $\cc\in L^r(\g;\H)$. Then
$T_t({\rm div}_\g\cc)$ is a function in $L^r(\g)$ for all $t>0$.
\end{lemma}
\begin{proof}
We use Proposition~\ref{pOU}(iii) to obtain
$$
|\langle T_t({\rm div}_\g\cc),\phi\rangle|= |\langle {\rm
div}_\g\cc,T_t\phi\rangle|\leq \int_E \|\cc\|_\H\|\nabla
T_t\phi\|_\H\,d\g\leq
C(q,t)\|\cc\|_{L^r(\g;\H)}\|\phi\|_{L^{r'}(\g)}
$$
for all $\phi\in\Cyl{E,\g}$, and we conclude.
\end{proof}

In the sequel we shall denote by $(\Lambda(p))^p$ the $p$-th
moment of the standard Gaussian in $\R$, i.e.
\begin{equation}\label{defLambdap}
\Lambda(p):=\biggl((2\pi)^{-1/2}\int_\R
|x|^pe^{-|x|^2/2}\,dx\biggr)^{1/p}.
\end{equation}

\begin{proposition}[Commutator estimate]\label{pcommu}
Let $\cc\in L^p(\g;\H)\cap LD^q(\g;\H)$ for some $p>1$, $1 \leq
q\leq 2$, with ${\rm div}_\g\cc\in L^q(\g)$. Let $r=\max\{p',q'\}$
and set
\begin{equation}\label{defre}
r^\e=r^\e(v,\cc):=e^\e \langle \cc,\nabla T_\e (v)\rangle-T_\e
({\rm div}_\g (v\cc)).
\end{equation}
Then, for $\e>0$ and $v\in L^r(\g)$ we have
\begin{equation}\label{maincommu}
\|r^\e\|_{L^1(\g)}\leq
\|v\|_{L^r(\g)}\biggl[\frac{\Lambda(p)\e}{\sqrt{1-e^{-2\e}}}\Vert
\cc\Vert_{L^p(\g;\H)}+\sqrt{2}\|{\rm div}_\g \cc\|_{L^q(\g)}+ 2\|
\|(\nabla\cc)^{\rm sym}\|_{HS}\|_{L^q(\g)}\biggr].
\end{equation}
Finally, $-r^\e\to v{\rm div}_\g\cc$ in $L^1(\g)$ as $\e\downarrow
0$.
\end{proposition}

\begin{proof}
The a-priori estimate \eqref{maincommu}, which is indeed the main
technical point of this paper, will be proved in the
Section~\ref{sect:finitedim} in finite-dimensional spaces. Here we
will just mention how the finite-dimensional approximation can be
performed.

Let us first assume that $v \in L^\infty$. Since $v\cc\in
L^p(\g;\H)$, the previous lemma ensures that $r^\e$ is a function.
Keeping $\cc$ fixed, we see that if $v_n\to v$ strongly in
$L^r(\g)$ then $r^\e(v_n,\cc)\to r^\e(v,\cc)$ in the duality with
$\Cyl{E,\g}$, and since the $L^1(\g)$ norm is lower semicontinuous
with respect to convergence in this duality, thanks to the density
of cylindrical functions we see that it suffices to prove
\eqref{maincommu} when $v$ is cylindrical. Keeping now
$v\in\Cyl{E,\g}$ fixed, we consider the vector fields
$$
\cc_N:=\sum_{i=1}^N\E_N\cc^i e_i.
$$
We observe that \eqref{defgdiv} gives ${\rm div}_\g\cc_N=\E_N({\rm
div}_\g\cc)$, while \eqref{explsymm} gives $(\n \cc_N)^{\rm
sym}=\E_N(\n \cc)^{\rm sym}$. Thus, by Jensen's inequality for
conditional expectations we obtain
$\Vert\cc_N\Vert_{L^p(\g;\H)}\leq\Vert\cc\Vert_{L^p(\g;\H)}$ and
$$
\int_E|{\rm div}_\g\cc_N|^q\,d\g\leq\int_E|{\rm
div}_\g\cc|^q\,d\g, \qquad \int_E\|(\nabla\cc_N)^{\rm
sym}\|_{HS}^q\,d\g\leq \int_E\|(\nabla\cc)^{\rm
sym}\|_{HS}^q\,d\g.
$$
Now, assuming that $v$ depends only on $\langle
e_1^*,x\rangle,\ldots,\langle e_M^*,x\rangle$, if we choose a
cylindrical test function $\phi$ depending only on $\langle
e_1^*,x\rangle,\ldots,\langle e_N^*,x\rangle$, with $N\geq M$
(with no loss of generality, because $v$ is fixed), we get
\begin{eqnarray*}
&&\int_E r^\e(v,\cc)\phi\,d\g=\int_E r^\e(v,\cc_N)\phi\,d\g\leq
\sup |\phi|\int_E
|r^\e(v,\cc_N)|\,d\g\\
&\leq& \sup|\phi|
\|v\|_{L^r(\g)}\biggl[\frac{\Lambda(p)\e}{\sqrt{1-e^{-2\e}}}\Vert
\cc_N\Vert_{L^p(\g;\H)}+\sqrt{2}\|{\rm div}_\g \cc_N\|_{L^q(\g)}+
2\| \|(\nabla\cc_N)^{\rm
sym}\|_{HS}\|_{L^q(\g)}\biggr]\\
&\leq& \sup|\phi|
\|v\|_{L^r(\g)}\biggl[\frac{\Lambda(p)\e}{\sqrt{1-e^{-2\e}}}\Vert
\cc\Vert_{L^p(\g;\H)}+\sqrt{2}\|{\rm div}_\g \cc\|_{L^q(\g)}+ 2\|
\|(\nabla\cc)^{\rm sym}\|_{HS}\|_{L^q(\g)}\biggr].
\end{eqnarray*}
This means that, once we know \eqref{maincommu} in
finite-dimensional spaces, we obtain that the same inequality
holds in all Wiener spaces for all $v\in L^\infty(\g)$. Finally,
to remove also this restriction on $v$, we consider a sequence
$(v_n)\subset L^\infty(\g)$ converging in $L^r(\g)$ to $v$ and we
notice that, because of \eqref{maincommu}, $r^\e(v_n,\cc)$ is a
Cauchy sequence in $L^1$ converging in the duality with
$\Cyl{E,\g}$ to $r^\e(v,\cc)$.

The strong convergence of $r^\e$ can be achieved by a density
argument. More precisely, if $q>1$ (so that $r<\i$), since
$r^\e(v,\cc)=r^\e(v-\phi,\cc)+r^\e(\phi,\cc)$, by
\eqref{maincommu} and the density of cylindrical functions in
$L^r(\g)$, we need only to consider the case when $v=\phi$ is
cylindrical. In this case
$$
r^\e=\langle \cc,T_\e\nabla\phi\rangle-T_\e(\phi{\rm
div}_\g\cc+\langle\cc,\nabla\phi\rangle)
$$
and its convergence to $-\phi{\rm div}_\g\cc$ is an obvious
consequence of the continuity properties of $T_\e$.

In the case $q=1$ (that is $r=\infty$), the approximation argument
is a bit more involved. Since we will never consider
$L^\infty$-regular flows, we give here just a sketch of the proof.
We argue as in \cite{lions3}: we write
$r^\e(v,\cc)=r^\e(v,\cc-\tilde\cc)+r^\e(v-\tilde v,\tilde
\cc)+r^\e(\tilde v,\tilde\cc)$, with $\tilde v$ and $\tilde \cc$
smooth and bounded with all their derivatives. Using
\eqref{maincommu} twice, we first choose $\tilde \cc$ so that
$r^\e(v,\cc-\tilde\cc)$ is small uniformly in $\e$, and then,
since now $\tilde \cc$ is smooth with bounded derivatives, it
suffices to choose $\tilde v$ close to $v$ in $L^s$ for some $s>1$
to make $r^\e(v-\tilde v,\tilde \cc)$ small. We can now conclude
as above.
\end{proof}

The following lemma is standard (both properties can be proved by a
smoothing argument; for the second one, see
\cite[Corollary~5.4.3]{Bogachev}):

\begin{lemma}[Chain rules]\label{lchain}
Let $\beta\in C^1(\R)$ with $\beta'$ bounded.
\begin{itemize}
\item[(i)] If $u,\,f\in L^1\bigl(I;L^1(\g)\bigr)$ satisfy
$\frac{d}{dt}u=f$ in the weak sense, then
$\frac{d}{dt}\beta(u)=\beta'(u)f$, still in the weak sense.
\item[(ii)] If $u\in W^{1,p}_\H(\g)$ then $\beta(u)\in
W^{1,p}_\H(\g)$ and $\nabla\beta(u)=\beta'(u)\nabla u$.
\end{itemize}
\end{lemma}

\begin{theorem}[Renormalization property]\label{main3}
Let $\bb:I\times E\to\H$ be satisfying the assumptions of the
uniqueness part of Theorem~\ref{main2}, with $I$ in place of
$(0,T)$. Then any solution $u_t$ of the continuity equation
\eqref{conti} in $L^\infty\bigl(I;L^r(\g)\bigr)$, with
$r=\max\{p',q'\}$, is renormalized.
\end{theorem}
\begin{proof} In the first step we prove the renormalized property
assuming that $u_t\in W^{1,r}_\H(\g)$ for a.e. $t$, and that both
$u_t$ and $\|\nabla u_t\|_\H$ belong to
$L^\infty\bigl(I;L^r(\g)\bigr)$. Under this assumption,
Remark~\ref{rprodu} gives that ${\rm div}_\g(\bb_t u_t)=u_t{\rm
div}_\g\bb_t+\langle\bb_t,\nabla u_t\rangle_\H$, therefore
$$
\frac{d}{dt}u_t=-u_t{\rm div}_\g\bb_t+\langle\bb_t,\nabla
u_t\rangle_\H\in L^1\bigl(I;L^1(\g)\bigr).
$$
Now, using Lemma~\ref{lchain} and Remark~\ref{rprodu} again, we
get
\begin{eqnarray*}
\frac{d}{dt}\beta(u_t)&=&
-\beta'(u_t)u_t{\rm div}_\g\bb_t-\beta'(u_t)\langle\bb_t,\nabla u_t\rangle_\H\\
&=& [\beta(u_t)-\beta'(u_t)u_t]{\rm div}_\g\bb_t-\beta(u_t){\rm
div}_\g\bb_t-
\langle\bb_t,\nabla\beta(u_t)\rangle_\H\\
&=& [\beta(u_t)-\beta'(u_t)u_t]{\rm div}_\g\bb_t-{\rm
div}_\g(\bb_t\beta(u_t)).
\end{eqnarray*}

Now we prove the renormalization property in the general case. Let
us define $u^\e_t:= e^{-\e} T_\e(u_t)$; since $T_\e$ is
self-adjoint in the sense of Proposition~\ref{pOU}(ii) and $T_\e$
maps cylindrical functions into cylindrical functions, the
continuity equation $\frac{d}{dt}u_t+{\rm div}_\g(\bb_tu_t)=0$
gives, still in the weak sense of duality with cylindrical
functions,
$$
\frac{d}{dt}u^\e_t+e^{-\e} T_\e [{\rm div}_\g (\bb_t u_t)]=0.
$$
Recalling the definition \eqref{defre}, we may write
$$
\frac{d}{dt}u^\e_t+{\rm div}_\g(\bb_tu^\e_t)=
e^{-\e}r^\e+u^\e_t{\rm div}_\g\bb_t.
$$
Denoting by $f^\e$ the right hand side, we know from
Proposition~\ref{pcommu} that $f^\e\to 0$ in
$L^1\bigl((0,T);L^1(\g))$. Taking into account that $u^\e_t$ and
$\|\nabla u^\e_t\|_\H$ belong to $L^\infty\bigl(I;L^r\g)\bigr)$
(by Proposition~\ref{pOU}(iii)), from the first step we obtain
$$
\frac{d}{dt}\beta(u^\e_t)+{\rm div}_\g(\bb_t\beta(u^\e_t))=
[\beta(u^\e_t)-u^\e_t\beta'(u^\e_t)]{\rm div}_\g\bb_t
+\beta'(u^\e_t)f^\e
$$
for all $\beta\in C^1(\R)$ with $\beta'(z)$ and
$z\beta'(z)-\beta(z)$ bounded. So, passing to the limit as
$\e\downarrow 0$ we obtain that $u_t$ is a renormalized solution.
\end{proof}

\noindent {\bf Proof of Theorem~\ref{main2}.} (Existence) It can be
obtained as a byproduct of the results in Section~\ref{sect:ODE}:
Theorem~\ref{teflows1} provides a generalized flow, i.e. a positive
finite measure $\eeta$ in the space of paths $\Omega(E)$, whose
marginals $(e_t)_\#\eeta$ at all times have a density uniformly
bounded in $L^r(\g)$, and $(e_0)_\#\eeta=\bar u\g$. Then, denoting
by $u_t$ the density of $(e_t)_\#\eeta$ with respect to $\g$,
Proposition~\ref{pindco} shows that $u_t$ solve the continuity
equation.

(Uniqueness) By the linearity of the equation, it suffices to show
that $\bar u=0$ implies $u_t\leq 0$ for all $t\in [0,T]$ for all
solutions $u\in L^\infty\bigl((0,T);L^r(\g)\bigr)$. We extend $u_t$
and $\bb_t$ to the interval $I:=(-1,T)$ by setting $u_t=\bar u$ and
$\bb_t=0$ for all $t\in (-1,0]$, and it is easy to check that this
extension preserves the validity of the continuity equation (still
in the weak form).

We choose, as a $C^1$ approximation of the positive part, the
functions $\beta_\e(z)$ equal to $\sqrt{z^2+\e^2}-\e$ for $z\geq
0$, and null for $z\leq 0$. Thanks to Theorem~\ref{main3}, we can
apply \eqref{conti1} with $\beta=\beta_\e$, with the test function
$\phi\equiv 1$, to obtain
$$
\frac{d}{dt}\int_E \beta_\e(u_t)\,d\g= \int_E
[\beta_\e(u_t)-u_t\beta_\e'(u_t)]{\rm div}_\g \bb_t\,d\g\leq
\e\int_E [{\rm div}_\g\bb_t]^-\,d\g,
$$
where we used the fact that $-\e\leq \beta_\e(z)-z\beta_\e'(z)\leq
0$. Letting $\e\downarrow 0$ we obtain that $\frac{d}{dt}\int_E
u_t^+\,d\g\leq 0$ in $(-1,T)$ in the sense of distributions. But
since $u_t=0$ for all $t\in (-1,0)$, we obtain $u_t^+=0$ for all
$t\in [0,T)$.

\section{Existence, uniqueness and stability of the flow}\label{sect:ODE}

In this section we discuss the problems of existence and uniqueness
of a flow associated to $\bb:[0,T]\times E\to\H$, and we discuss its
main properties.

\subsection{Existence of a generalized $\bb$-flow}

It will be useful, in order to establish our first existence
result, a definition of flow more general than
Definition~\ref{dflow1}. In the sequel we shall denote by
$\Omega(E)$ the space of continuous maps from $[0,T]$ to $E$,
endowed with the sup norm. Since $E$ is separable, $\Omega(E)$ is
complete and separable. We shall denote by
$$
e_t:\Omega(E)\to E,\qquad e_t(\omega):=\omega(t)
$$
the evaluation maps at time $t\in [0,T]$.

If $1\leq\alpha\leq\infty$, we shall also denote by
$AC^\alpha(E)\subset\Omega(E)$ the subspace of functions $\omega$
satisfying
\begin{equation}\label{finx}
\omega(t)=\omega(0)+\int_0^t g(s)\,ds\qquad\forall t\in [0,T]
\end{equation}
for some $g\in L^\alpha\bigl((0,T);E\bigr)$. The function $g$,
that we shall denote by $\dot\omega$, is uniquely determined up to
negligible sets by \eqref{finx}: indeed, if $\bar t$ is a Lebesgue
point of $g$ then $\langle e^*,g(\bar t)\rangle$ coincides with
the derivative at $t=\bar t$ of the real-valued absolutely
continuous function $t\mapsto\langle e^*,\omega(t)\rangle$, for
all $e^*\in E^*$.

\begin{definition}[Generalized $\bb$-flows and $L^r$-regularity]\label{dflow2}
If $\bb:[0,T]\times E\to E$, we say that a probability measure
$\eeta$ in $\Omega(E)$ is a flow associated to $\bb$ if:
\begin{itemize}
\item[(i)] $\eeta$ is concentrated on maps $\omega\in AC^1(E)$ satisfying
the ODE $\dot\omega=\bb(t,\omega)$ in the integral sense, namely
\begin{equation}\label{bochner1}
\omega(t)=\omega(0)+\int_0^t \bb_\tau(\omega(\tau))\,d\tau \qquad
\forall t\in [0,T];
\end{equation}
\item[(ii)] $(e_0)_\#\eeta=\g$.
\end{itemize}
If in addition there exists $1\leq r\leq\infty$ such that, for all
$t\in [0,T]$, the image measures $(e_t)_\#\eeta$ are absolutely
continuous with respect to $\g$ with a density in $L^r(\g)$, then we
say that the flow is $L^r$-regular.
\end{definition}

\begin{remark}[Invariance of $\bb$-flows]\label{invariaflow}
{\rm Assume that $\eeta$ is a generalized $L^1$-regular $\bb$-flow
and $\tilde{\bb}$ is a modification of $\bb$, i.e., for a.e. $t\in
(0,T)$ the set $N_t:=\{\bb_t\neq\tilde{\bb}_t\}$ is $\g$-negligible.
Then, because of $L^1$-regularity we know that, for a.e. $t\in
(0,T)$, $\omega(t)\notin N_t$ $\eeta$-almost surely. By Fubini's
theorem, we obtain that, for $\eeta$-a.e. $\omega$, the set of times
$t$ such that $\omega(t)\in N_t$ is negligible in $(0,T)$. As a
consequence $\eeta$ is a $\tilde{\bb}$-flow as well. }\end{remark}

\begin{remark}[Martingale solutions of ODEs]
{\rm We remark that the notion of generalized flow coincides with
the Stroock-Varadhan's notion of martingale solutions for
stochastic differential equations in the particular case when
there is no noise (so that the stochastic differential equation
reduces to an ordinary differential equations), see for instance
\cite{strvar} and \cite[Lemma 3.8]{figalli}.}\end{remark}
From now on, we shall adopt the convention $\|v\|_\H=+\infty$ for
$v\in E\setminus\H$.

\begin{proposition}[Compactness]\label{ptight}
Let $K\subset E$ be a compact set, $C\geq 0$, $\alpha\in
(1,\infty)$ and let $\mathcal F\subset AC^\alpha(E)$ be the family
defined by:
$$
\mathcal F:=\left\{\omega\in AC^\alpha(E)\,:\,\omega(0)\in
K,\,\,\, \int_0^T\|\dot\omega\|_\H^\alpha\,dt\leq C\right\}.
$$
Then $\mathcal F$ is compact in $\Omega(E)$.
\end{proposition}
\begin{proof} Let us fix an integer $h$, and split $[0,T]$ in the
$h$ equal intervals $I_i:=[iT/h,(i+1)T/h]$, $i=0,\ldots,h-1$. We
consider the family ${\mathcal F}_h$ obtained by replacing each
curve $\omega(t)$ in $\mathcal F$ with the continuous ``piecewise
affine'' curve $\omega_h$ coinciding with $\omega$ at the
endpoints of the intervals $I_i$ and with constant derivative,
equal to $\frac{T}{h}\int_{I_i}\dot\omega(t)\,dt$, in all
intervals $(iT/h,(i+1)T/h)$. We will check that each family
$\mathcal F_h$ is relatively compact, and that
$\sup|\omega-\omega_h|\to 0$ as $h\to\infty$, uniformly with
respect to $\omega\in\mathcal F$. These two facts obviously imply,
by a diagonal argument, the relative compactness of $\mathcal F$.

The family $\mathcal F_h$ is easily seen to be relatively compact:
indeed, the initial points of the curve lie in the compact set
$K$, and since $\{\int_{I_0}\dot\omega(t)\,dt\}_{\omega\in\mathcal
F}$ is uniformly bounded in $\H$, the compactness of the embedding
of $\H$ in $E$ shows that also the family of points
$\{\omega_h(T/h)\}_{\omega\in\mathcal F}$ is relatively compact;
continuing in this way, we prove that all families of points
$\{\omega_h(iT/h)\}_{\omega\in\mathcal F}$, $i=0,\ldots,h-1$, and
therefore the family ${\mathcal F}_h$, are relatively compact.

Fix $\omega\in\mathcal F$; denoting by $L$ the norm of the
embedding of $\H$ in $E$, we have
$$
\|\omega(t)-\omega_h(t)\|\leq\int_{iT/h}^t
\|\dot\omega(\tau)-\dot\omega_h(\tau)\|\,d\tau \leq 2L
\int_{iT/h}^{t}\|\dot\omega(\tau)\|_\H\,d\tau\leq
2LC^{1/\alpha}\biggl(\frac{T}{h}\biggr)^{1-1/\alpha}
$$
for all $t\in [iT/h,(i+1)T/h]$. This proves the uniform
convergence of $\omega_h$ to $\omega$ as $h\to\infty$, as $\omega$
varies in $\mathcal F$.

Finally, we have to check that $\mathcal F$ is closed. The
stability of the condition $\omega(0)\in K$ under uniform
convergence is obvious. The stability of the second condition can
be easily obtained thanks to the reflexivity of the space
$L^\alpha\bigl((0,T);\H\bigr)$.
\end{proof}

\begin{theorem}[Existence of $L^r$-regular generalized $\bb$-flows]\label{teflows1}
Let $\bb:[0,T]\times E\to\H$ be satisfying the assumptions of the
existence part of Theorem~\ref{main2}. Then there exists a
generalized $\bb$-flow $\eeta$, $L^r$-regular for all $r\in
[1,c/T]$. In addition, the density $u_t$ of $(e_t)_\#\eeta$ with
respect to $\g$ satisfies
\begin{equation}\label{buondi1}
\int (u_t)^r\,d\g\leq\biggl\Vert\int\exp\bigl(Tr [{\rm
div}_\g\bb_t]^-\bigr)\,d\g\biggr\Vert_{L^\infty(0,T)} \qquad \forall
t \in [0,T].
\end{equation}
\end{theorem}
\begin{proof} {\bf Step 1.} (finite-dimensional approximation)
Let $\bb_N:[0,T]\times E\to\H_N$ be defined by
$\sum_{i=1}^N\bb_N^i e_i$, where
$$
\bb_N^i(t,\cdot):=\E_N\bb^i_t,\qquad 1\leq i\leq N,\,\,t\in [0,T].
$$
Arguing as in the proof of Proposition~\ref{pcommu}, we have the
estimates
\begin{equation}\label{ccases}
\int_0^T \biggl(\int_E \|(\bb_N)_t\|_\H^p
\,d\g(x)\biggr)^{1/p}\,dt\leq \int_0^T\biggl(\int_E
\|\bb_t\|_\H^p\,d\g(x)\biggr)^{1/p}\,dt,
\end{equation}
\begin{equation}\label{ccases2}
\biggl\|\int_E\exp\bigl(c[{\rm div}_\g
(\bb_N)_t]^-\bigr)\,d\g(x)\biggr\|_{L^\i(0,T)}\leq \biggl\|\int_E
\exp\bigl(c[{\rm
div}_\g\bb_t]^-\bigr)\,d\g(x)\biggr\|_{L^\i(0,T)}.
\end{equation}

By applying Theorem~\ref{texflow} to the finite-dimensional fields
$\tilde{\bb}_N$ given by the restriction of $\bb_N$ to
$[0,T]\times\H_N$, we obtain a generalized flow $\ssigma_N$ in
$\H_N$ (i.e. a positive finite measure in $\Omega(H_N)$) associated
to $\tilde{\bb}_N$. Using the inclusion map $i_N$ of $\H_N$ in $\H$
we obtain a generalized flow $\eeta_N:=(i_N)_\#\ssigma_N$ associated
to $\bb_N$. In addition, \eqref{ccases2} and the finite-dimensional
estimate \eqref{finalstim} give
\begin{equation}\label{ccases1}
\sup_{t\in [0,T]}\sup_{N\geq 1}\int_E (u^N_t)^r\,d\g\leq \|\bar
u\|_{L^\infty(\g)}^r \biggl\Vert\int_E\exp\bigl(Tr [{\rm
div}_\g\bb_t]^-\bigr)\,d\g\biggr\Vert_{L^\infty(0,T)},
\end{equation}
with $u^N_t$ equal to the density of $(e_t)_\#\eeta_N$ with
respect to $\g$.

\noindent {\bf Step 2.} (Tightness and limit flow $\eeta$).
 We call coercive a functional $\Psi$ if its sublevel sets
 $\{\Psi\leq C\}$ are compact.
 Since $(\E_N\bar u\g)$ is a tight family of measures,
  by Prokhorov theorem we can find (see
 for instance \cite{strvar}) a coercive functional $\Phi_1:E\to
 [0,+\infty)$ such that $\sup_N\int_E \Phi_1\E_N\bar u\,d\g<\infty$.
 We choose $\alpha\in (1,p)$ such that $(p/\alpha)'\leq c/T$ (this is possible
 because we are assuming that $p'T<c$) and consider the functional
 \begin{equation}\label{deffhi}
 \Phi(\omega):=
 \begin{cases}
 \Phi_1(\omega(0))+\int_0^T\|\dot\omega(t)\|_\H^\alpha\,dt &\text{if
 $\omega\in AC^p(E)$;}
 \\
 +\infty &\text{if $\omega\in\Omega(E)\setminus AC^\alpha(E)$.}
 \end{cases}
\end{equation}
Thanks to Proposition~\ref{ptight} and the coercivity of $\Phi_1$,
$\Phi$ is a coercive functional in $\Omega(E)$. Since
\begin{eqnarray*}
\int_{\Omega(E)}\Phi(\omega)\,d\eeta_N(\omega)&=&
\int_E\Phi_1(x)\E_N\bar u(x)\,d\g(x)+\int_0^T\int_{\Omega(E)}
\|(\bb_N)_t(\omega(t))\|_\H^\alpha\,d\eeta_N(\omega)\,dt\\&=&
\int_E\Phi_1(x)\E_N\bar u(x)\,d\g(x)+\int_0^T\int_E
\|(\bb_N)_t(x)\|_\H^\alpha u_t^N(x)\,d\g(x)\,dt
\end{eqnarray*}
we can apply H\"older inequality with the exponents $p/\alpha$ and
$(p/\alpha)'$, \eqref{ccases}, \eqref{ccases2} and \eqref{ccases1}
to obtain that $\int\Phi\,d\eeta_N$ is uniformly bounded. So, we can
apply again Prokhorov theorem to obtain that $(\eeta_N)$ is tight in
$\Omega(E)$. Therefore we can find a positive finite measure $\eeta$
in $\Omega(E)$ and a family of integers $N_i\to\infty$ such that
$\eeta_{N_i}\to\eeta$ weakly, in the duality with
$C_b\bigl(\Omega(E)\bigr)$. In the sequel, to simplify our notation,
we shall assume that convergence occurs as $N\to\infty$. Obviously,
because of \eqref{ccases1}, $\eeta$ is $L^r$-regular and, more
precisely, \eqref{buondi1} holds.

\noindent {\bf Step 3.} ($\eeta$ is a $\bb$-flow).
 It suffices to show that
\begin{equation}\label{valleggi}
\int_{\Omega(E)} 1\land
\|\omega(t)-\omega(0)-\int_0^t\bb_s(\omega(s))\,ds\|\,d\eeta=0
\end{equation}
for any $t\in [0,T]$. The technical difficulty is the integrand in
\eqref{valleggi}, due to the lack of regularity of $\bb_t$, is not
continuous in $\Omega(E)$; the truncation with the constant 1 is
used to have a bounded integrand. To this aim, we prove first that
\begin{equation}\label{valleggi1}
\int_{\Omega(E)} 1\land
\|\omega(t)-\omega(0)-\int_0^t\cc_s(\omega(s))\,ds\|\,d\eeta\leq
\int_0^T\int_E\|\bb_s(x)-\cc_s(x)\|u_s(x)\,d\g(x)\, ds
\end{equation}
for any bounded continuous function $\cc$. Then, choosing a
sequence $(\cc_n)$ converging to $\bb$ in
$L^1\bigl((0,T);L^p(\g;E)\bigr)$, and noticing that
$$
\int_{\Omega(E)}\int_0^T\|\bb_s(\omega(s))-(\cc_n)_s(\omega(s))\|
\,ds\, d\eeta
=\int_0^T\int_E\|\bb_s(x)-(\cc_n)_s(x)\|u_s(x)\,d\g(x)\,ds
\rightarrow 0,
$$
we can pass to the limit in \eqref{valleggi1} with $\cc=\cc_n$ to
obtain \eqref{valleggi}.

It remains to show \eqref{valleggi1}. This is a limiting argument
based on the fact that \eqref{valleggi} holds for $\bb_N$,
$\eeta_N$:
\begin{eqnarray*}
&&\int_{\Omega(E)} 1\land
\|\omega(t)-\omega(0)-\int_0^t\cc_s(\omega(s))\,ds\|\,d\eeta
\\&=&\lim_{N\to\infty} \int_{\Omega(E)}1\land
\|\omega(t)-\omega(0)-\int_0^t\cc_s(\omega(s))\,ds\|\,d\eeta_N\\&=&
\lim_{N\to\infty}\int_{\Omega(E)} 1\land\|\int_0^t
(\bb_N)_s(\omega(s))-\cc_s(\omega(s))\,ds\|\,d\eeta_N
\\&\leq&\limsup_{N\to\infty}
\int_0^T\int_E \|(\bb_N)_s(x)-\cc_s(x)\|u^N_s(x)\,d\g(x)\,ds
=\int_0^T\int_E \|\bb_s(x)-\cc_s(x)\|u_s(x)\,d\g(x)\,ds.
\end{eqnarray*}
In order to obtain the last equality we added and subtracted
$\|\bb_s-\cc_s\|u^N_s$, and we used the strong convergence of
$\bb_N$ to $\bb$ in $L^1\bigl((0,T);L^p(\g;E)\bigr)$ and the
weak$^*$ convergence of $u^N_s$ to $u_s$ in
$L^\i\bigl((0,T);L^{p'}(\g;E)\bigr)$.
\end{proof}

\subsection{Uniqueness of the $\bb$-flow}

The following lemma provides a simple characterization of Dirac
masses (i.e. measures concentrated at a single point), for measures
in $\Omega(E)$ and for families of measures in $E$.

\begin{lemma}\label{simple}
Let $\ssigma$ be a positive finite measure in $\Omega(E)$. Then
$\ssigma$ is a Dirac mass if and only if $(e_t)_\#\ssigma$ is a
Dirac mass for all $t\in\Q\cap [0,T]$. \\
A Borel family $\{\nu_x\}_{x\in E}$ of positive finite measures in
$E$ (i.e. $x\mapsto\nu_x(A)$ is Borel in $E$ for all $A\subset E$
Borel) is made, for $\g$-a.e. $x$, by Dirac masses if and only if
\begin{equation}\label{orthomu}
\nu_x(A_1)\nu_x(A_2)=0\quad\text{$\g$-a.e. in $E$, for all
disjoint Borel sets $A_1,\,A_2\subset E$.}
\end{equation}
\end{lemma}
\begin{proof}
The first statement is a direct consequence of the fact that all
elements of $\Omega(E)$ are continuous maps, which are uniquely
determined on $\Q\cap [0,T]$. In order to prove the second
statement, let us fix an integer $k\geq 1$ and a countable
partition $(A_i)$ of Borel sets with ${\rm diam}(A_i)\leq 1/k$
(its existence is ensured by the separability of $E$). By
\eqref{orthomu} we obtain a $\g$-negligible Borel set $N_k$
satisfying $\nu_x(A_i)\nu_x(A_j)=0$ for all $x\in E\setminus N_k$.
As a consequence, the support of each of the measures $\nu_x$, as
$x$ varies in $E\setminus N_k$, is contained in the closure of one
of the sets $A_i$, which has diameter less than $1/k$. It follows
that $\nu_x$ is a Dirac mass for all $x\in
E\setminus\bigcup_kN_k$.
\end{proof}

\begin{theorem}[Uniqueness of $L^r$-regular generalized
$\bb$-flows]\label{tgflows2} Let $\bb:[0,T]\times E\to\H$ be
satisfying the assumptions of the uniqueness part of
Theorem~\ref{main2}, let $r=\max\{p',q'\}$ and assume that $c\geq
rT$. Let $\eeta$ be a $L^r$-regular generalized $\bb$-flow. Then:
\begin{itemize}
\item[(i)] for $\g$-a.e. $x\in E$, the measures
$\E(\eeta|\omega(0)=x)$ are Dirac masses in $\Omega(E)$, and
setting
\begin{equation}\label{defX}
\E(\eeta|\omega(0)=x)=\delta_{\sxX(\cdot,x)},\qquad
\XX(\cdot,x)\in\Omega(E),
\end{equation}
the map $\XX(t,x)$ is a $L^r$-regular $\bb$-flow, according to
Definition~\ref{dflow1}. \item[(ii)] Any other $L^r$-regular
generalized $\bb$-flow coincides with $\eeta$. In particular $\XX$
is the unique $L^r$-regular $\bb$-flow.
\end{itemize}
\end{theorem}
\begin{proof} (i) We set $\eeta_x:=\E(\eeta|\omega(0)=x)$.
Taking into account the first statement in Lemma~\ref{simple}, it
suffices to show that, for $\bar t\in\Q\cap [0,T]$ fixed, the
measures $\nu_x:=\E((e_{\bar t})_\#\eeta|\omega(0)=x)=(e_{\bar
t})_\#\eeta_x$ are Dirac masses for $\g$-a.e. $x\in E$. Still
using Lemma~\ref{simple}, we will check the validity of
\eqref{orthomu}. Since $\nu_x=\delta_x$ when $\bar t=0$, we shall
assume that $\bar t>0$.

Let us argue by contradiction, assuming the existence of a Borel
set $L\subset E$ with $\g(L)>0$ and disjoint Borel sets
$A_1,\,A_2\subset E$ such that both $\nu_x(A_1)$ and $\nu_x(A_2)$
are positive for $x\in L$. We will get a contradiction with
Theorem~\ref{main2}, building two distinct solutions of the
continuity equation with the same initial condition $\bar u\in
L^\infty(\g)$. With no loss of generality, possibly passing to a
smaller set $L$ still with positive $\g$-measure, we can assume
that the quotient $\beta(x):=\nu_x(A_1)/\nu_x(A_2)$ is uniformly
bounded in $L$. Let $\Omega_i\subset\Omega(E)$ be the set of
trajectories $\omega$ which belong to $A_i$ at time $\bar t$;
obviously $\Omega_1\cap\Omega_2=\emptyset$ and we can define
positive finite measures $\eeta_i$ in $\Omega(E)$ by
$$
\eeta_1:=\int_L\chi_{\Omega_1}\eeta_x\,d\g(x),
\qquad\eeta_2:=\int_L \beta(x)\chi_{\Omega_2}\eeta_x\,d\g(x).
$$
By Proposition~\ref{pindco}, both $\eeta_1$ and $\eeta_2$ induce,
via the identity $u^i_t\g=(e_t)_\#\eeta_i$, a solution to the
continuity equation which is uniformly bounded (just by comparison
with the one induced by $\eeta$) in $L^r(\g)$. Moreover, both
solutions start from the same initial condition $\bar
u(x)=\nu_x(A_1)\chi_L(x)$. On the other hand, by the definition of
$\Omega_i$, $u^1_{\bar t}\g$ is concentrated in $A_1$ while
$u^2_{\bar t}\g$ is concentrated in $A_2$, therefore $u^1_{\bar
t}\neq u^2_{\bar t}$. So, uniqueness of solutions to the
continuity equation is violated.

\noindent (ii) If $\ssigma$ is any other $L^r$-regular generalized
$\bb$-flow, we may apply statement (i) to the flows $\ssigma$, to
obtain that for $\g$-a.e. $x$ also the measures
$\E(\ssigma|\omega(0)=x)$ are Dirac masses; but since the property
of being a generalized flow is stable under convex combinations,
also the measures
$$
\frac{1}{2}\E(\eeta|\omega(0)=x)+
\frac{1}{2}\E(\ssigma|\omega(0)=x)=
\E\bigl(\frac{\eeta+\ssigma}{2}|\omega(0)=x\bigr)
$$
must be Dirac masses for $\g$-a.e. $x$. This can happen only if
$\E(\eeta|\omega(0)=x)=\E(\ssigma|\omega(0)=x)$ for $\g$-a.e. $x$.
\end{proof}

The connection between solutions to the ODE $\dot \XX=\bb_t(\XX)$
and the continuity equation is classical: in the next proposition
we present it under natural regularity assumptions in this
setting.

\begin{proposition}\label{pindco}
Let $\eeta$ be a positive finite measure in $\Omega(E)$
satisfying:
\begin{itemize}
\item[(a)] $\eeta$ is concentrated on paths $\omega\in AC^1(E)$
such that $\omega(t)=\omega(0)+\int_0^t\bb_s(\omega(s))\,ds$ for
all $t\in [0,T]$; \item[(b)]
$\int_0^T\int_{\Omega(E)}\|\dot\omega(t)\|_\H\,d\eeta(\omega)\,dt<\infty$.
\end{itemize}
Then the measures $\mu_t:=(e_t)_\#\eeta$ satisfy
$\frac{d}{dt}\mu_t+{\rm div}_\g(\bb_t\mu_t)=0$ in $(0,T)\times E$
in the weak sense.
\end{proposition}
\begin{proof} Let $\phi(x)=\psi(\langle e_1^*,x\rangle,\ldots,\langle e_N^*,x\rangle)$
be cylindrical. By (a) and Fubini's theorem, for a.e. $t$ the
following property holds: the maps $\langle
e_i^*,\omega(t)\rangle$, $1\leq i\leq N$, are differentiable at
$t$, with derivative equal to $\langle
e_i^*,\bb_t(\omega(t))\rangle$, for $\eeta$-a.e. $\omega$. Taking
\eqref{svitla} into account, for a.e. $t$ we have
\begin{eqnarray*}
\frac{d}{dt}\int_E \phi\,d\mu_t&=& \frac{d}{dt}\int_{\Omega(E)}
\psi(\langle e_1^*,\omega(t)\rangle,\ldots,\langle
e_N^*,\omega(t)\rangle)\,d\eeta\\&=&
\sum_{i=1}^N\int_{\Omega(E)}\frac{\partial\psi}{\partial z_i}
(\langle e_1^*,\omega(t)\rangle,\ldots,\langle e_N^*,\omega(t)\rangle)\langle e_i^*,\dot \omega(t)\rangle\,d\eeta\\
&=& \sum_{i=1}^N\int_{\Omega(E)}\frac{\partial\psi}{\partial z_i}
(\langle e_1^*,\omega(t)\rangle,\ldots,\langle
e_N^*,\omega(t)\rangle)\langle
e_i,\bb_t(\omega(t))\rangle_\H\,d\eeta\\&=& \int_E
\langle\nabla\phi,\bb_t\rangle_\H\,d\mu_t.
\end{eqnarray*}
In the previous identity we used, to pass to the limit under the
integral sign, the property
$$
\lim_{h\to 0}\langle e_i^*,\frac{\omega(t+h)-\omega(t)}{h}\rangle=
\langle e_i^*,\dot\omega(t)\rangle\quad\text{in $L^1(\eeta)$, for
$1\leq i\leq N$},
$$
whose validity for a.e. $t$ is justified by assumption (b). The same
assumption also guaranteees (see for instance \cite[\S3]{cetraro}
for a detailed proof) that $t\mapsto \int_E\phi\,d\mu_t$ is
absolutely continuous, so its pointwise a.e. derivative coincides
with the distributional derivative.
\end{proof}

\subsection{Stability of the $\bb$-flow and semigroup property}

The methods we used to show existence and uniqueness of the flow
also yield stability of the flow with respect to approximations (not
necessarily finite-dimensional ones) of the vector field. In the
proof we shall use the following simple lemma (see for instance
Lemma~22 of \cite{cetraro} for a proof), where we use the notation
${\rm id}\times f$ for the map $x\mapsto (x,f(x))$.

\begin{lemma}[Convergence in law and in probability]\label{rollo}
Let $F$ be a metric space and let $f_n,\,f:E\to F$ be Borel maps.
Then $f_n\to f$ in $\g$-probability if and only if ${\rm id}\times
f_n\to {\rm id}\times f$ in law.
\end{lemma}

\begin{theorem}[Stability of $L^r$-regular $\bb$-flows]\label{main4}
Let $p,\,q>1$, $r=\max\{p',q'\}$ and let $\bb_n,\,\bb:(0,T)\times
E\to\H$ be satisfying:
\begin{itemize}
\item[(i)] $\bb_n\to\bb$ in $L^1\bigl((0,T);L^p(\g;\H)\bigr)$;
\item[(ii)] for a.e. $t\in (0,T)$ we have $(\bb_n)_t,\,\bb_t\in
LD^q_\H(\g;\H)$ with
\begin{equation}\label{ldsymbis}
\sup_{n\in\N}\int_0^T\biggl(\int_E \|(\nabla(\bb_n)_t)^{\rm
sym}(x)\|^q_{HS}\,d\g(x)\biggr)^{1/q}\, dt<\infty
\end{equation}
and ${\rm div}_\g(\bb_n)_t$ and ${\rm div}_\g\bb_t$ belong to
$L^1\bigl((0,T);L^q(\g)\bigr)$; \item[(iii)] $\exp (c[{\rm div}_\g
(\bb_n)_t]^-)$ are uniformly bounded in
$L^\infty\bigl((0,T);L^1(\g)\bigr)$ for some $c\geq T r$.
\end{itemize}
Then, denoting by $\XX_n$ (resp. $\XX$) the unique $L^r$ regular
$\bb_n$-flows (resp. $\bb$-flow) we have
\begin{equation}\label{inmisura}
\lim_{n\to\infty}\int_E\sup_{[0,T]}\|\XX_n(\cdot,x)-\XX(\cdot,x)\|\,d\g(x)=0.
\end{equation}
\end{theorem}
\begin{proof}
Let us denote the generalized $\bb_n$-flows $\eeta_n$ induced by
$\XX_n$, namely the law under $\g$ of $x\mapsto\XX_n(\cdot,x)$.
The uniform estimates (iii), together with the boundedness of
$\|\bb_n\|_\H$ in $L^1\bigl((0,T);L^p(\g)\bigr)$ imply, in view of
\eqref{buondi1},
\begin{equation}\label{buondi}
\sup_{n\in\N}\int
(u_t^n)^r\,d\g\leq\sup_{n\in\N}\biggl\Vert\int\exp\bigl(Tr [{\rm
div}_\g\bb^n_t]^-\bigr)\,d\g\biggr\Vert_{L^\infty(0,T)}<\infty
\qquad \forall t \in [0,T],
\end{equation}
where $u_t^n$ is the density of $(e_t)_\#\eeta_n=X(t,\cdot)_\#\g$
with respect to $\g$. In addition, by the same argument used in Step
2 of the proof of Theorem~\ref{teflows1} we have
$$
\sup_{n\in\N}\int_{\Omega(E)}\Phi(\omega)\,d\eeta_n(\omega)<\infty,
$$
where $\Phi$ is defined as in \eqref{deffhi}, with $\alpha\in (1,p)$
and $\Phi_1:E\to [0,\infty)$ $\gamma$-integrable and coercive.

This estimate implies the tightness of $(\eeta_n)$. If $\eeta$ is a
limit point, in the duality with $C_b(\Omega(E))$, of $\eeta_n$, the
same argument used in Step 3 of the proof of Theorem~\ref{teflows1}
gives that $\eeta$ is a generalized $\bb$-flow. In addition, the
uniform estimates \eqref{buondi} imply that $\eeta$ is
$L^r$-regular. As a consequence we can apply Theorem~\ref{tgflows2}
to obtain that $\eeta$ is the law of the $\Omega(E)$-valued map
$x\mapsto\XX(\cdot,x)$, and more precisely that
$\E(\eeta|\omega(0)=x)=\delta_{\sxX(\cdot,x)}$ for $\g$-a.e. $x$.
Therefore, by the uniqueness of $\XX$, the whole sequence
$(\eeta_n)$ converges to $\eeta$ and $\XX_n$ converge in law to
$\XX$.

In order to obtain that $x\mapsto\XX_n(\cdot,x)$ converge in
$\g$-probability to $x\mapsto\XX(\cdot,x)$ we use Lemma~\ref{rollo}
with $F=\Omega(E)$, so we have to show that ${\rm
id}\times\XX_n(\cdot,x)$ converge in law to ${\rm
id}\times\XX(\cdot,x)$. For all $\psi\in
C_b\bigl(E\times\Omega(E)\bigr)$ we have
$$
\int_E\psi(x,\XX_n(\cdot,x))\,d\g(x)= \int_{\Omega(E)}
\psi(e_0(\omega),\omega)\,d\eeta_n\to
\int_{\Omega(E)}\psi(e_0(\omega),\omega)\,d\eeta= \int_E
\psi(x,\XX(\cdot,x))\,d\g(x),
$$
and this proves the convergence in law.

Finally, by adding and subtracting $x$, we can prove
\eqref{inmisura} provided we show that
$\sup_{[0,T]}|\XX(\cdot,x)-x|\in L^1(\g)$ and
$\sup_{[0,T]}|\XX_n(\cdot,x)-x|$ are equi-integrable in $L^1(\g)$.
We prove the second property only, because the proof of the first
one is analogous. Starting from the integral formulation of the ODE,
Jensen's inequality gives $\sup_{[0,T]}|\XX_n(\cdot,x)-x|^\alpha\leq
T^{\alpha-1}\int_0^T\|\bb_\tau(\XX_n(\tau,x))\|\,d\tau$ and by
integrating both sides with respect to $\g$, Fubini's theorem gives
$$
\int_E\sup_{[0,T]}|\XX_n(\cdot,x)-x|^\alpha\,d\g(x)\leq
T^{\alpha-1}\int_E\int_0^T\int_E\|\bb_\tau\|^\alpha
u^n_\tau\,d\g\,d\tau.
$$
Choosing $\alpha>1$ such that $(p/\alpha)'\leq c/T$ (this is
possible because we are assuming that $c>p'T$) and applying the
H\"older inequality with the exponents $p/\alpha$ and $(p/\alpha)'$
we obtain that $\sup_{[0,T]}|\XX_n(\cdot,x)-x|$ are equibounded in
$L^\alpha(\g)$.
\end{proof}

Under the same assumptions of Theorem~\ref{tgflows2}, for all $s\in
[0,T]$ also a unique $L^r$-regular flow $\XX^s:[s,T]\times E\to E$
exists, characterized by the properties that
$\tau\mapsto\XX^s(\tau,x)$ is an absolutely continuous map in
$[s,T]$ satisfying
\begin{equation}\label{XsX}
\XX^s(t,x)=x+\int_s^t\bb_\tau\bigl(\XX^s(\tau,x)\bigr)\,d\tau\qquad
\forall t\in [s,T]
\end{equation}
for $\g$-a.e. $x\in E$, and the regularity condition
$\XX^s(\tau,\cdot)_\#\g=f_\tau\g$, with $f_\tau\in L^r(\g)$, for all
$\tau\in [s,T]$. This family of flow maps satisfies the semigroup
property:

\begin{proposition}[Semigroup property]\label{psemi}
Under the same assumptions of Theorem~\ref{tgflows2}, the unique
$L^r$-regular flows $\XX^s$ starting at time $s$ satisfy the
semigroup property
\begin{equation}\label{siver1}
\XX^s\left(t,\XX^r(s,x)\right)=\XX^r(t,x)
\qquad\text{for $\g$-a.e. $x\in E$, for all $0\leq r\leq s\leq t\leq
T$.}
\end{equation}
\end{proposition}
\begin{proof} Let $r,\,s,\,t$ be fixed. By combining the finite-dimensional
projection argument of Step 1 of the proof of
Theorem~\ref{teflows1}, with the smoothing argument used in Step 2
of the proof of Theorem~\ref{texflow} we can find a family of
vector fields $\bb_n$ converging to $\bb$ in
$L^1\bigl((0,T);L^p(\g;\H)\bigr)$ and satisfying the uniform
bounds of Theorem~\ref{main4}, whose (classical) flows $\XX_n$
satisfy the semigroup property (see \eqref{siver3})
\begin{equation}\label{siver}
\XX_n^s\left(t,\XX_n^r(s,x)\right)=\XX_n^r(t,x)
\qquad\text{for $\g$-a.e. $x\in E$, for all $0\leq r\leq s\leq t\leq
T$.}
\end{equation}
We will pass to the limit in \eqref{siver}, to obtain
\eqref{siver1}. To this aim, notice that \eqref{inmisura} of
Theorem~\ref{main4} immediately provides the convergence in
$L^1(\g)$ of the right hand sides, so that we need just to show
convergence in $\g$-measure of the left hand sides. Notice first
that the convergence in $\g$-measure of $\XX_n^r(s,\cdot)$ to
$\XX^r(s,\cdot)$ implies the convergence in $\g$-measure of
$\psi(\XX_n^r(s,\cdot))$ to $\psi(\XX^r(s,\cdot))$ for any Borel
function $\psi:E\to\R$ (this is a simple consequence of the fact
that, by Lusin's theorem, we can find a nondecreasing sequence of
compact sets $K_n\subset E$ such that $\psi\vert_{K_n}$ is uniformly
continuous and $\g(E\setminus K_n)\downarrow 0$, and of the fact
that the laws of $\XX_x^r(s,\cdot)$ are uniformly bounded in
$L^r(\g)$), so that choosing $\psi(z):=\XX^s(t,z)$, and adding and
subtracting $\XX^s(t,\XX_n(s,x))$, the convergence in $\g$-measure
of the right hand sides of \eqref{siver} to
$\XX^s\left(t,\XX^r(s,x)\right)$ follows by the convergence in
$\g$-measure to $0$ of
$$
\XX_n^s\left(t,\XX_n^r(s,x)\right)-\XX^s\left(t,\XX_n^r(s,x)\right).
$$
Denoting by $\rho_n$ the density of the law of $\XX_n^r(s,\cdot)$,
we have
$$
\int_E1 \wedge \|
\XX_n^s\left(t,\XX_n^r(s,x)\right)-\XX^s\left(t,\XX_n^r(s,x)\right)\|\,d\g(x)=
\int_E 1\wedge \|\XX_n^s(t,y)-\XX^s(t,y)\|\rho_n(y)\,d\g(y),
$$
and the right hand side tends to $0$ thanks to \eqref{inmisura} and
to the equi-integrability of $(\rho_n)$.
\end{proof}

The semigroup property allows also to construct a unique family of
flows $\XX^s:[s,T]\times E\times E$ even in the case when the
assumption \eqref{peeters} is replaced by
$$
\exp(c[{\rm div}_\g \bb_t]^-)\in
L^\infty\bigl((0,T);L^1(\g)\bigr)\qquad\text{for \emph{some} $c>0$.}
$$
The idea is to compose the flows defined on sufficiently short
intervals, with length $T'$ satisfying $c>rT'$. It is easy to check
that this family of flow maps is uniquely determined by the
semigroup property \eqref{siver1} and by the \emph{local} regularity
property
$$
\text{$\XX^s(t,\cdot)_\#\g\ll\g$ with a density in $L^r(\g)$ for all
$t\in [s,\min\{s+T',T\}]$, $s\in [0,T]$.}
$$
Globally in time, the only property retained is
$\XX^s(t,\cdot)_\#\g\ll\g$ for all $t\in [s,T]$.

\subsection{Convergence of finite-dimensional
flows}\label{finitedim}

Assume that we are given vector fields
$\bb_N:[0,T]\times\R^N\to\R^N$ satisfying, for some $p,\,q>1$ the
assumptions (i), (ii), (iii) of Theorem~\ref{main1} (with
$E=\H=\R^N$) relative to the standard Gaussian $\g_N$ in $\R^N$,
with norms uniformly bounded by constants independent of $N$. Let
us assume that $\bb_N$ is a consistent family, namely the
conditional expectation of the projection of $(\bb_{N+1})_t$ on
$\R^N$, given $x^1,\ldots,x^N$, is $(\bb_N)_t$. Let
$\XX_N:[0,T]\times\R^N\to\R^N$ be the associated $\bb_N$-flows.

In this section we briefly illustrate how the stability results of
this paper can be used to prove the convergence of $\XX_N$ and to
characterize their limit.

To this aim, let us denote by $\gamma_p$ the product of standard
Gaussians in the countable product $\R^\infty$, and notice that
the consistency assumption provides us with a unique vector field
$\bb:[0,T]\times\R^\infty\to\R^\infty$ such that, denoting by
$\E_N$ the conditional expectation with respect to
$x^1,\ldots,x^N$ and by $\pi_N:\R^\infty\to\R^N$ the canonical
projections, the identities $\E_N\pi_N\bb_t=(\bb_N)_t$ hold. In
order to recover a Wiener space we fix a sequence
$(\lambda_i)\in\ell^2$ and define
$$
E:=\left\{(x^i):\ \sum_{i=1}^\infty
\lambda_i^2(x^i)^2<\infty\right\}.
$$
The space $E$ can be endowed with the canonical scalar product,
and obviously $\gamma_p(E)=1$, so that $\bb$ can be also viewed as
a vector field in $E$ and the induced measure $\gamma$ in $E$ is
Gaussian. According to Remark~\ref{rhilbert}, its Cameron-Martin
space $\H$ can be identified with $\ell^2$. Then, we can apply the
stability Theorem~\ref{main4} (viewing, with a slight abuse,
$\bb_N$ as vector fields in $E$ and, consequently, their flows
$\XX_N$ as flows in $E$ which leave $x^{N+1},x^{N+2},\ldots$
fixed) to obtain that $\XX_N$ converge to the flow $\XX$ relative
to $\bb$ in $L^1(\g;E)$. It follows that
\begin{equation}\label{gross}
\lim_{N\to\infty}\int_{\R^\infty}\sqrt{\sum_{i=1}^\infty\lambda_i^2
|\XX_N^i(t,x)-\XX^i(t,x)|^2}\,d\g_p(x)=0 \qquad\forall t\in
[0,T],\,\,\forall (\lambda_i)\in\ell^2.
\end{equation}
Finally, notice that also $\XX$ could be defined without an
explicit mention to $E$, working in $(\R^\infty,\gamma_p)$ in
place of $(E,\gamma)$. According to this viewpoint, $E$ plays just
the role of an auxiliary space, and deliberately we wrote
\eqref{gross} without an explicit mention to it.

\section{An extension to non $\H$-valued vector fields}
\label{sect:nonCM}

In \cite{Peters}, \cite{bogachev}, the authors consider the
following equation:
\begin{equation}
\label{ODErotation} \XX(t,x)=\tilde Q_t x + \int_0^t
Q_{t-s}\bb_s(\XX(s,x))\,ds.
\end{equation}
Here $(Q_t)_{t \in \R}$ is a strongly continuous group of orthogonal
operator on $\H$, and $\tilde Q_t :E \to E$ denotes the measurable
linear extension of $Q_t$ to $E$ (which always exists and preserves
the measure $\g$, see for instance \cite{kusuoka}). Observe that,
thanks to the Duhamel formula, \eqref{ODErotation} formally
corresponds to the equation
$$
\dot\XX(t,x)=L\XX(t,x) + \bb_t(\XX(t,x)),
$$
where $L$ denotes the generator of the group (i.e. $\dot Q_t=LQ_t$).

The definition of $L^r$-regular flow can be extended in the
obvious way to \eqref{ODErotation}. Let us now see how our results
allow to prove existence and uniqueness of $L^r$-regular flows
under the assumptions of Theorem~\ref{main1} (observe that this
forces in particular $r>1$).

Let $\XX(t,x)$ be a solution of \eqref{ODErotation}, and define
$\YY(t,x):=\tilde Q_{-t}\XX(t,x)$. Then we have
\begin{eqnarray*}
\YY(t,x)&=&x + \int_0^t Q_{-s}\bb_s(\XX(s,x))\,ds\\
&=&x + \int_0^t Q_{-s}\bb_s(\tilde Q_s\YY(s,x))\,ds.
\end{eqnarray*}
Therefore $\YY$ is a flow associated to the vector field
$\cc_t(x):=Q_{-t}\bb_t(\tilde Q_t x)$. Moreover $\YY$ is still a
$L^r$-regular flow. Indeed, if $u_t\in L^r(\g)$ denotes the density
of the law of $\XX(t,\cdot)$, then, for all $\phi\in\Cyl{E,\g}$, we
have
\begin{eqnarray*}
\int \phi(\YY(t,x))\,d\g(x)&=&\int \phi(\tilde
Q_{-t}\XX(t,x))\,d\g(x)=\int \phi(\tilde Q_{-t}x)u_t(x)\,d\g(x)\\
&\leq &\|u_t\|_{L^r(\g)} \|\phi \circ \tilde Q_t\|_{L^{r'}(\g)}
=\|u_t\|_{L^r(\g)} \|\phi \|_{L^{r'}(\g)}.
\end{eqnarray*}
Since $r>1$, this implies that $\YY$ is $L^r$-regular. On the other
hand we remark that, using the same argument, one obtains that, if
$\YY$ is a $L^r$-regular flow associated to $\cc$, then
$\XX(t,x):=\tilde{Q}_t\YY(t,x)$ is a $L^r$-regular flow for
\eqref{ODErotation}.

We have therefore shown that there is a one-to-one correspondence
between $L^r$-regular flows for \eqref{ODErotation} and
$L^r$-regular flows associated to $\cc$. To conclude the existence
and uniqueness of $L^r$-regular flows for \eqref{ODErotation}, it
suffices to observe that, thanks to the orthogonality of $Q_t$ and
the measure-preserving property of $\tilde Q_t$, if $\bb$ satisfies
all the assumptions in Theorem~\ref{main1}, then so does $\cc$
thanks to the identities $\|\cc_t(x)\|_\H=\|\bb_t(\tilde{Q}_t
x)\|_\H$, $\|(\n \cc_t)^{\rm sym}(x)\|_{HS}=\|(\n\bb_t)^{\rm
sym}(\tilde{Q}_t x)\|_{HS}$, and $\div_\g
\cc_t(x)=\div_\g\bb_t(\tilde{Q}_tx)$.

Indeed, let us check the formula for the symmetric part of the
derivative, the proof of the one concerning the divergence being
similar and even simpler. Let $h=j(e^*)\in\H$ and notice that $Q_th
=j(f^*)$, where $\langle f^*,y\rangle=\langle
e^*,\tilde{Q}_{-t}(y)\rangle$. Using Remark~\ref{rinvaria} and the
fact that $\phi\mapsto\phi\circ\tilde{Q}_t$ maps $\Cyl{E,\gamma}$
into $C^1_b(E,\g)$, for $\phi\in\Cyl{E,\g}$ we get
\begin{eqnarray*}
&&\int_E\langle\cc_t,h\rangle_\H\partial_h\phi\,d\g= \int_E\langle
\bb_t(\tilde{Q}_tx),Q_th\rangle_\H\partial_h\phi(x)\,d\g(x)\\&=&
\int_E\langle
\bb_t(y),Q_th\rangle_\H(\partial_h\phi)\circ\tilde{Q}_{-t}(y)\,d\g(y)=
\int_E\langle
\bb_t(y),Q_th\rangle_\H\partial_{Q_th}(\phi\circ\tilde{Q}_{-t})(y)\,d\g(y)\\
&=&-\int_E\partial_{Q_th}\langle\bb_t,Q_th\rangle_\H\phi\circ\tilde{Q}_{-t}\,d\g(y)
+\int_E\langle \bb_t(y),Q_th\rangle_\H\phi\circ\tilde{Q}_{-t}\langle
f^*,y\rangle\,d\g(y)\\
&=&
-\int_E[\partial_{Q_th}\langle\bb_t,Q_th\rangle_\H]\circ\tilde{Q}_t\phi\,d\g(x)
+\int_E\langle\cc_t(x),h\rangle_\H\phi\langle e^*,x\rangle\,d\g(x).
\end{eqnarray*}
This proves that
$\partial_h\langle\cc_t,h\rangle_\H=
\partial_{Q_th}\langle\bb_t,Q_th\rangle_\H\circ\tilde{Q}_t$,
and using the fact that $Q_t$ maps orthonormal bases of $\H$ in
orthonormal bases of $\H$ we get $\|(\n\cc_t)^{\rm
sym}\|_{HS}=\|(\n\bb_t)^{\rm sym}\|_{HS}\circ\tilde{Q}_t$.

\section{Finite-dimensional estimates}
\label{sect:finitedim}

This section is devoted to the proof of the crucial a-priori
bounds \eqref{expbound} and \eqref{maincommu} in
finite-dimensional Wiener spaces. So, we shall assume that
$E=\H=\R^N$ and, only in this section, denote by $x\cdot y$ the
scalar product in $\R^N$, and by $|x|$ the Euclidean norm
(corresponding to the norm of the Cameron-Martin space). Also,
only in this section we shall denote by $\g$ the standard Gaussian
in $\R^N$, product of $N$ standard Gaussians in $\R$, and by
$\int$ integrals on the whole of $\R^N$. The sums $\sum_i$ (resp.
$\sum_{i,j}$) will always be understood with $i$ (resp. $i$ and
$j$) running from $1$ to $N$.

\subsection{Upper bounds on the flow density}\label{sect:flow}

In this subsection we show the existence part of
Theorem~\ref{main2} in finite-dimensional Wiener spaces
$E=\H=\R^N$.

\begin{theorem}\label{texflow}
Let $\bb:(0,T)\times\R^N\to\R^N$ be satisfying the assumptions of
the existence part of Theorem~\ref{main2}. Then, for any $r \in
[1,c/T]$ there exists a generalized $L^r$-regular $\bb$-flow
$\eeta$. Its density $u_t$ satisfies also
\begin{equation}\label{finalstim}
\int (u_t)^r\,d\g\leq\biggl\Vert\int\exp\bigl(Tr [{\rm
div}_\g\bb_t]^-\bigr)\,d\g\biggr\Vert_{L^\infty(0,T)} \qquad \forall
t \in [0,T].
\end{equation}
\end{theorem}
\begin{proof} {\bf Step 1.} Here we consider first the case when $\bb_t$ are
smooth, with $\int_0^T\|\nabla\bb_t\|_{L^\infty(B)}\,dt$ finite for
all bounded open sets $B\subset\R^N$. Under this assumption, for all
$x\in\R^N$ the unique solution $\XX(\cdot,x)$ to the ODE $\dot
\XX(t,x)=\bb_t(\XX(t,x))$, with the initial condition $\XX(0,x)=x$,
is defined until some maximal time $\tau(x)\in (0,T]$. Obviously, by
the maximality of $\tau(x)$, if
$$
\limsup_{t\uparrow\tau(x)}|\XX(t,x)|<+\infty
$$
then $\tau(x)=T$ and the solution is continuous in $[0,T]$.

Let us fix $s \in [0,T)$. We denote $E_s$ the set $\{\tau>s\}$ and
notice that standard stability results for ODE's with a locally
Lipschitz vector field ensure that $E_s$ is open and that
$x\mapsto\XX(t,x)$ is smooth in $E_s$ for $t\in [0,s]$. Furthermore,
from the identity
$\dot\nabla_x\XX(t,x)=\nabla\bb_t(\XX(t,x))\nabla_x\XX(t,x)$,
obtained by spatial differentiation of the ODE (see \cite{cetraro}
for details), one obtains
\begin{equation}\label{binet}
\dot J\XX(t,x)={\rm div}\bb_t(\XX(t,x))J\XX(t,x)\qquad x\in
E_s,\,\,t\in [0,s],
\end{equation}
where $J\XX(t,x)$ is the determinant of $\nabla_x\XX(t,x)$.

We first compute a pointwise expression for the measure
$\XX(t,\cdot)_\#(\chi_{E_s}\g)$ for $t\in [0,s]$. By the change of
variables formula, the density $\rho_t^s$ of
$\XX(t,\cdot)_\#(\chi_{E_s}\g)$ with respect to $\Leb{N}$ is linked
to the initial density $\bar\rho^s$ by
$$
\rho_t^s(\XX(t,x))=\frac{\bar\rho^s(x)}{J\XX(t,x)},
$$
where $\bar\rho^s(y):=\chi_{E_s}(y) e^{-|y|^2/2}$. Denoting by
$u^s_t$ the density of $\XX(t,\cdot)_\#(\chi_{E_s}\g)$ with respect
to $\g$, we get
\begin{equation}\label{basicODE}
u^s_t\bigl(\XX(t,x)\bigr)=\frac{\bar\rho^s(x)}{J\XX(t,x)}e^{|\sxX(t,x)|^2/2}.
\end{equation}
So, taking the identity \eqref{binet} into account, we obtain
$$
\frac{d}{dt}u^s_t\bigl(\XX(t,x)\bigr)=-{\rm
div}_\g\bb_t\bigl(\XX(t,x)\bigr)
\frac{\bar\rho^s(x)}{J\XX(t,x)}e^{|\sxX(t,x)|^2/2}= -{\rm
div}_\g\bb_t\bigl(\XX(t,x)\bigr)u^s_t\bigl(\XX(t,x)\bigr).
$$
By integrating the ODE, for $t\in [0,s]$ we get
\begin{eqnarray*}
u^s_t\bigl(\XX(t,x)\bigr)&=&\chi_{E_s}(x)\exp\biggl(-\int_0^t{\rm
div}_\g\bb_\tau\bigl(\XX(\tau,x)\bigr)\,d\tau\biggr)\\
&\leq&\chi_{E_s}(x) \exp\biggl(\int_0^t[{\rm
div}_\g\bb_\tau\bigl(\XX(\tau,x)\bigr)]^-\,d\tau\biggr).
\end{eqnarray*} We can now estimate $\Vert u^s_t\Vert_{L^r(\g)}$ as
follows:
\begin{eqnarray*}
\int (u_t^s)^r\,d\g&=&\int (u_t^s)^{r-1} u^s_t\,d\g\leq
\int\exp\biggl((r-1)\int_0^t[{\rm
div}_\g\bb_\tau\bigl(\XX(\tau,x)]^-\bigr)\,d\tau\biggr)\chi_{E_s}(x)\,d\g(x)\\
&\leq& \int\frac{1}{t}\int_0^t\exp\bigl(t(r-1)[{\rm
div}_\g\bb_\tau\bigl(\XX(\tau,x)\bigr)]^-\bigr)\,d\tau\chi_{E_s}(x)\,d\g(x)\\&=&
\frac{1}{t}\int_0^t\int\exp\bigl(t(r-1)[{\rm
div}_\g\bb_\tau\bigl(\XX(\tau,x)\bigr)]^-\bigr)\chi_{E_s}(x)\,d\g(x)\,d\tau\\
&\leq& \frac{1}{t}\int_0^t\int\exp\bigl(T(r-1)[{\rm
div}_\g\bb_\tau(y)]^-\bigr)u_\tau^s(y)\,d\g(y)\,d\tau.
\end{eqnarray*}
Now, set $\Lambda(t):=\int_0^t\Vert
u^s_\tau\Vert^r_{L^r(\g)}\,d\tau$ and apply the H\"older
inequality to get
\begin{eqnarray}\label{perugia}
\Lambda'(t)&\leq&\frac{1}{t}\biggl(\int_0^t\int\exp\bigl( Tr [{\rm
div}_\g\bb_\tau(y)]^-\bigr)\,d\g(y)\,d\tau\biggr)^{1/r'}\Lambda^{1/r}(t)
\\&\leq& K t^{1/r'-1}
\Lambda^{1/r}(t)=Kt^{-1/r}\L^{1/r}(t),\nonumber
\end{eqnarray}
with $K:=\Vert\int\exp\bigl(T r [{\rm
div}_\g\bb_t]^-\bigr)\,d\g\Vert_{L^\infty(0,T)}^{1/r'}$. An
integration of this differential inequality yields $\Lambda(t)\leq
K^{r'} t$, which inserted into \eqref{perugia} gives
\begin{equation}\label{finalstim5}
\int (u_t^s)^r\,d\g\leq\biggl\Vert\int\exp\bigl(Tr [{\rm
div}_\g\bb_t]^-\bigr)\,d\g\biggr\Vert_{L^\infty(0,T)}
\qquad\forall t\in [0,s], \, \forall s \in [0,T).
\end{equation}

Now, let us prove that the flow is globally defined in $[0,T]$ for
$\g$-a.e. $x$: we have indeed
\begin{eqnarray*}
\int\sup_{[0,\tau(x))}|\XX(t,x)-x|\,d\g(x)&\leq&
\int\int_0^{\tau(x)}|\bb_t(\XX(t,x))|\,dt\,d\g(x)=
\int_0^T\int_{E_t}|\bb_t(\XX(t,x))|\,d\g(x)\,dt\\&=&
\int_0^T\int|\bb_t|u^t_t\,d\g\,dt.
\end{eqnarray*}
Using \eqref{finalstim5} with $s=t$, we obtain that
$\int\sup\limits_{[0,\tau(x))}|\XX(t,x)-x|\,d\g(x)$ is finite, so
that $\tau(x)=T$ and $\XX(\cdot,x)$ is continuous up to $t=T$ for
$\g$-a.e. $x$. Letting $s\uparrow T$ in \eqref{finalstim5} we obtain
\eqref{finalstim}.

Denoting as in \eqref{XsX} by $\XX^s$ the flow starting at time $s$,
we also notice (this is useful in the proof, by approximation, of
the semigroup property in Proposition~\ref{psemi}) that the
pointwise uniqueness of the flow implies the semigroup property
\begin{equation}\label{siver3}
\XX^s\left(t,\XX^r(s,x)\right)=\XX^r(t,x) \qquad\text{for all $0\leq
r\leq s\leq t\leq T$}
\end{equation}
for all $x$ where $\XX^r(\cdot,x)$ is globally defined in $[r,T]$.

\noindent {\bf Step 2.} In this step we remove the regularity
assumptions made on $\bb$, considering the vector fields $\bb_\e$
defined by $\bb_\e^i(t,\cdot):=T_\e\bb^i_t$. It is immediate to
check that the fields $\bb_\e$ satisfy the regularity assumptions
made in Step 1, so the existence of a $L^r$-regular $\bb_\e$-flow
$\eeta_\e$ satisfying
\begin{equation}\label{finalstim2}
\int (u_t^\e)^r\,d\g\leq\biggl\Vert\int\exp\bigl(Tr [{\rm div}_\g
(\bb_\e)_t]^-\bigr)\,d\g\biggr\Vert_{L^\infty(0,T)}
\end{equation}
is ensured by Step 1. In \eqref{finalstim2} the functions $u_t^\e$
are, as usual, the densities of $(e_t)_\#\eeta_\e$ with respect to
$\g$. Now, since ${\rm div}_\g((\bb_\e)_t)=e^{-\e} T_\e({\rm
div}_\g\bb_t)$, we may apply Jensen's inequality to get
\begin{equation}\label{finalstim3}
\int (u_t^\e)^r\,d\g\leq\biggl\Vert\int\exp\bigl( e^{-\e} Tr [{\rm
div}_\g \bb_t]^-\bigr)\,d\g\biggr\Vert_{L^\infty(0,T)}.
\end{equation}
Since
$$
\int_0^T\biggl(\int\|\bb_\e(t,x)\|_\H^p\,d\g\biggr)^{1/p}\,dt\leq
\int_0^T\biggl(\int\|\bb(t,x)\|_\H^p\,d\g\biggr)^{1/p}\,dt,
$$
the same tightness argument used in the proof of
Theorem~\ref{teflows1} to pass from finitely many to infinitely
many dimensions provides us with a $\bb$-flow $\eeta$ satisfying
\eqref{finalstim}: any weak limit point $\eeta$ of $\eeta_\e$ as
$\e\downarrow 0$.
\end{proof}

\subsection{Commutator estimate}\label{sect:commu}

This subsection is entirely devoted to the proof of the commutator
estimate \eqref{maincommu} in finite-dimensional Wiener spaces.

We will often use the ``Gaussian rotations''
\begin{equation}\label{Gausschange}
(x,y)\mapsto
(z,w):=\bigl(e^{-\e}x+\sqrt{1-e^{-2\e}}y,-\sqrt{1-e^{-2\e}}x+e^{-\e}y\bigr),
\end{equation}
mapping the product measure $\g(dx)\times\g(dy)$ into
$\g(dz)\times\g(dw)$. Indeed, the transformations above preserve
the Lebesgue measure in $\R^N\times\R^N$ (being their Jacobian
identically equal to 1) and $|x|^2+|y|^2=|z|^2+|w|^2$.

We now state two elementary Gaussian estimates. The first one
\begin{equation}\label{cancellaGauss1}
\left(\int |l\cdot w|^p\,d\g(w)\right)^{1/p}= |l|\left(\int
|w_1|^p\,d\g(w)\right)^{1/p}=\Lambda(p)|l|\qquad\forall l\in\R^N,
\end{equation}
with $\Lambda$ depending only on $p$, is a simple consequence of
the rotation invariance of $\g$.

\begin{lemma}
Let $A:\R^N\to\R^N$ be a linear map and $c\in\R$. Then, if $q\leq
2$, we have
\begin{equation}\label{cancellaGauss2}
\biggl(\int\bigl|\langle Aw,w\rangle
-c\bigr|^q\,d\g(w)\biggr)^{1/q}\leq \sqrt{2}\|A^{\rm
sym}\|_{HS}+|{\rm tr\,}A-c|.
\end{equation}
\end{lemma}
\begin{proof}
Obviously we can assume that $A$ is symmetric. By rotation
invariance, we can also assume that $A$ is diagonal, and denote by
$\lambda_1,\ldots,\lambda_N$ its eigenvalues. We have then
\begin{eqnarray*}
\int\bigl|\sum_i\lambda_i(w^i)^2-c\bigr|^2\,d\g(w)&=& \int
\Bigl[\sum_{ij}\lambda_i\lambda_j (w^i)^2(w^j)^2- 2c\sum_i
\lambda_i(w^i)^2+c^2\Bigr]\,d\g(w)\\
&=&
3\sum_i\lambda_i^2+\sum_{i\neq j}\lambda_i\lambda_j -2c\sum_i\lambda_i+c^2\\
&=&
2\sum_i\lambda_i^2+\sum_{ij}\lambda_i\lambda_j-2c\sum_i\lambda_i+c^2\\
&=& 2\sum_i\lambda_i^2+\bigl(\sum_i\lambda_i-c\bigr)^2.
\end{eqnarray*}
If $q=2$ we take the square roots of both sides and we conclude;
if $q\leq 2$ we apply the H\"older inequality.
\end{proof}

Henceforth, a vector field $\cc\in L^p(\g;\R^N)\cap
LD^q_\H(\g;\R^N)$ and a function $v\in L^r(\g)$ will be fixed,
with $r=\max\{p',q'\}$ and $p>1$, $1\leq q \leq 2$. Our goal is to
prove the estimate
\begin{equation}\label{maincommu1}
\|r^\e\|_{L^1(\g)}\leq
\|v\|_{L^r(\g)}\biggl[\frac{\Lambda(p)\e}{\sqrt{1-e^{-2\e}}}\Vert
\cc\Vert_{L^p(\g;\R^N)}+ 2^{1/q'}\|{\rm div}_\g \cc\|_{L^q(\g)}+
2^{1/q'}\sqrt{2}\| \|(\nabla\cc)^{\rm sym}\|_{HS}
\|_{L^q(\g)}\biggr],
\end{equation}
where
\begin{equation}\label{hilbo}
r^\e:=e^\e\cc\cdot\nabla v_\e-T_\e({\rm div}_\g(v\cc)).
\end{equation}
Since $2^{1/q'}\leq\sqrt{2}$, this yields the finite-dimensional
version of \eqref{maincommu}.

In this setup the Ornstein-Uhlenbeck operator $v_\e:=T_\e v$ takes
the explicit form
$$
v_\e(x):=\int v(e^{-\e}x + \sqrt{1-e^{-2\e}}y)\,d\g(y) =\int
v(z)\rho_\e(x,z)\,d\g(z)
$$
with
\begin{eqnarray*}
\rho_\e(x,z)&:=&\frac{1}{(1-e^{-2\e})^{N/2}}\exp(-\frac{|e^{-\e}x-z|^2}{2(1-e^{-2\e})})\exp(\frac{|z|^2}{2})
\\&=&\frac{1}{(1-e^{-2\e})^{N/2}}\exp(-\frac{|e^{-\e}x|^2
- 2 \e^{-\e}x \cdot z + |e^{-\e}z|^2}{2(1-e^{-2\e})}).
\end{eqnarray*}
This implies that
\begin{eqnarray} \label{nablave} \n v_\e(x)&=&\int
v(z)\n_x\rho_\e(x,z)\,d\g(z)=-e^{-\e}\int
\frac{e^{-\e}x-z}{1-e^{-2\e}}f(z)\rho_\e(x,z)\,d\g(z)\nonumber\\
&=&e^{-\e}\int v(e^{-\e}x + \sqrt{1-e^{-2\e}}y)
\frac{y}{\sqrt{1-e^{-2\e}}}\,d\g(y).
\end{eqnarray}
Let us look for a more explicit expression of the commutator in
\eqref{hilbo}. To this aim, we show first that $T_\e({\rm
div}_\g(v\cc))$ is a function, and
\begin{equation}\label{divmoll}
T_\e(\div_\g(v\cc ))(x)=\int (v\cc)(e^{-\e}x + \sqrt{1-e^{-2\e}}y)
\cdot\frac{y}{\sqrt{1-e^{-2\e}}}\,d\g(y) - T_\e(z\cdot v\cc)(x).
\end{equation}
If $\cc$ and $v$ are smooth, this is immediate to check: indeed,
thanks to \eqref{expdiv}, we need only to show that
$$
T_\e({\rm div\,}(v\cc ))(x)=\int (v\cc)(e^{-\e}x +
\sqrt{1-e^{-2\e}}y) \cdot\frac{y}{\sqrt{1-e^{-2\e}}}\,d\g(y).
$$
The latter is a direct consequence of \eqref{nablave} (with $v$
replaced by $v\cc^i$) and of the relation $\partial_i T_\e
(v\cc^i)=e^{-\e} T_\e(\partial_i(v\cc^i))$. If $v$ and $\cc$ are
not smooth, we argue by approximation.

Therefore, taking \eqref{nablave} and \eqref{divmoll} into
account, we have that $r^\e(x)$ is given by
\begin{align*}
&\int v(e^{-\e}x + \sqrt{1-e^{-2\e}}y)\frac{\cc(x)-\cc(e^{-\e}x +
\sqrt{1-e^{-2\e}}y)}{\sqrt{1-e^{-2\e}}}
\cdot y \,d\g(y)\\
&+\int v(e^{-\e}x + \sqrt{1-e^{-2\e}}y)\cc(e^{-\e}x +
\sqrt{1-e^{-2\e}}y)\cdot(e^{-\e}x +
\sqrt{1-e^{-2\e}}y)\,d\g(y)\\
&= \int \frac{v(e^{-\e}x +
\sqrt{1-e^{-2\e}}y)}{\sqrt{1-e^{-2\e}}}\left\{\cc(x)\cdot y -
\cc(e^{-\e}x+\sqrt{1-e^{-2\e}}y)\cdot(e^{-2\e}y-
e^{-\e}\sqrt{1-e^{-2\e}}x)\right\}\,d\g(y).
\end{align*}
Now, using the abbreviations
$\alpha_\e(x,y):=v(e^{-\e}x+\sqrt{1-e^{-2\e}}y)$,
$\beta_\e:=\e/\sqrt{1-e^{-2\e}}$, we interpolate and write
$-r^\e(x)$ as
\begin{align}\label{splittingre}
&\frac{1}{\sqrt{1-e^{-2\e}}}\int\alpha_\e(x,y)\frac{d}{dt}\int_0^1
\cc(e^{-t\e}x+\sqrt{1-e^{-2\e t}}y)\cdot (e^{-2t\e}
y-e^{-t\e}\sqrt{1-e^{-2t\e}}x)\,dt \,d\g(y)\nonumber \\
&= \beta_\e\int\alpha_\e(x,y)\\
&\int_0^1 \biggl[\sum_{ij}\left(\p_j \cc^i(e^{-t\e}x +
\sqrt{1-e^{-2t\e}}y)[e^{-t\e}\sqrt{1-e^{-2t\e}}x^i
-e^{-2t\e}y^i][e^{-t\e}x^j -
\frac{e^{-2t\e}}{\sqrt{1-e^{-2t\e}}}y^j]\right)\nonumber \\
&+\sum_i\left(\cc^i(e^{-t\e}x +
\sqrt{1-e^{-2t\e}}y)[(e^{-t\e}\sqrt{1-e^{-2t\e}}-
\frac{e^{-3t\e}}{\sqrt{1-e^{-2t\e}}})x^i
-2e^{-2t\e}y^i]\right)\biggr]\,dt\,d\g(y)\nonumber \\
&=:\beta_\e\int \alpha_\e(x,y) (A_\e(x,y)+B_\e(x,y))\,d\g(y),
\end{align}
where, adding and subtracting
$$
\sum_i\cc^i(e^{-t\e}x+\sqrt{1-e^{-2t\e}}y)\frac{e^{-2t\e}}{\sqrt{1-e^{-2t\e}}}
(e^{-t\e}x^i+\sqrt{1-e^{-2t\e}}y^i),
$$
we have set
\begin{eqnarray*}
A_\e(x,y):=\int_0^1\biggl(\sum_{ij}&\p_j \cc^i(e^{-t\e}x +
\sqrt{1-e^{-2t\e}}y)[e^{-t\e}\sqrt{1-e^{-2t\e}}x^i
-e^{-2t\e}y^i][e^{-t\e}x^j -\displaystyle{
\frac{e^{-2t\e}}{\sqrt{1-e^{-2t\e}}}}y^j]\nonumber\\
&{\displaystyle
-\sum_i\cc^i(e^{-t\e}x+\sqrt{1-e^{-2t\e}}y)\frac{e^{-2t\e}}{\sqrt{1-e^{-2t\e}}}
(e^{-t\e}x^i+\sqrt{1-e^{-2t\e}}y^i)\biggr)\,dt},
\end{eqnarray*}
$$
B_\e(x,y):=\int_0^1\sum_i\left(\cc^i(e^{-t\e}x +
\sqrt{1-e^{-2t\e}}y)e^{-t\e}
[\sqrt{1-e^{-2t\e}}x^i-e^{-t\e}y^i]\right)\,dt.
$$
Let us estimate $\beta_\e\int\int |\alpha_\e B_\e|\,d\g d\g$
first: the change of variables \eqref{Gausschange} and Fubini's
theorem give
$$
\beta_\e\int\int|\alpha_\e B_\e|\,d\g(x)\,
d\g(y)\leq\beta_\e\int_0^1 e^{-\e t}\int\int |v(z)|\bigl|
\sum_i\cc^i(z)w^i\bigr|\,d\g(z)\,d\g(w)\,dt.
$$
Using \eqref{cancellaGauss1} with $f=\cc(z)$, we get
\begin{equation}\label{finaleBe}
\beta_\e\int\int|\alpha_\e B_\e|\,d\g(x)\, d\g(y)\leq
\beta_\e\int\int |v(z)|\biggl| \sum_i\cc^i(z)w^i\biggr|\,d\g(z)\,d\g(w)\\
\leq\beta_\e\Lambda(p)\Vert \cc\Vert_{L^p(\g;\R^N)}\Vert
v\Vert_{L^{p'}(\g)}.
\end{equation}

Now, we estimate $\beta_\e\int\int |\a_\e A_\e|\,d\g\,d\g$; again,
we use the change of variables \eqref{Gausschange} to write
$$
e^{-t\e}\sqrt{1-e^{-2t\e}}x^i-e^{-2t\e}y^i= -e^{-t\e}w^i,\qquad
e^{-t\e}x^j-\frac{e^{-2t\e}}{\sqrt{1-e^{-2t\e}}}y^j=
-\frac{e^{-t\e}}{\sqrt{1-e^{-2t\e}}}w^j.
$$
Therefore we get
\begin{align*}
&\beta_\e\int\int|\alpha_\e A_\e|\,d\g(x)\,
d\g(y)\\
&\leq \beta_\e\int_0^1\int\int |v(z)|\biggl|\sum_{ij}\p_j \cc^i(z)
\frac{e^{-2t\e}}{\sqrt{1-e^{-2t\e}}}w^iw^j
-\sum_i\cc^i(z)\frac{e^{-2t\e}}{\sqrt{1-e^{-2t\e}}}z^i
\biggr|\,d\g(z)\,d\g(w)\,dt\\
&=\int\int
|v(z)|\biggl|\sum_{ij}\partial_j\cc^i(z)w^iw^j-\sum_i\cc^i(z)z^i\biggr|
\,d\g(z)\,d\g(w),
\end{align*}
where we used the identity
$$
\int_0^1\frac{e^{-2t\e}}{\sqrt{1-e^{-2t\e}}}\,dt
=\frac{\sqrt{1-e^{-2\e}}}{\e}=\beta_\e^{-1}.
$$
Eventually we use \eqref{cancellaGauss2} with $A=\nabla\cc(z)$ and
$c=\cc(z)\cdot z$ to obtain
\begin{eqnarray}\label{finaleAe}
\beta_\e\int\int|\alpha_\e A_\e|\,d\g(x)\,d\g(y)&\leq& \Vert
v\Vert_{L^{q'}(\g)}\biggl(\int\int
\bigl|\sum_{ij}\partial_j\cc^i(z)w^iw^j-\sum_i\cc^i(z)z^i\bigr|^q\,
d\g(w)\,d\g(z)\biggr)^{1/q}\nonumber \\
&\leq& 2^{1-1/q}\Vert v\Vert_{L^{q'}(\g)}\biggl(\int
\sqrt{2}^q\|\|(\nabla\cc)^{\rm sym}\|_{HS}\|^q+|{\rm
div}_\g\cc|^q\,d\g(z)\biggr)^{1/q}\nonumber
\\
&\leq& 2^{1-1/q}\|v\|_{L^{q'}(\g)}\biggl(\sqrt{2}\|
\|(\nabla\cc)^{\rm sym}\|_{HS}\|_{L^q(\g)}+ \|{\rm
div}_\g\cc\|_{L^q(\g)}\biggr).
\end{eqnarray}
Combining \eqref{splittingre}, \eqref{finaleBe} and
\eqref{finaleAe}, we have proved \eqref{maincommu1}.

\end{document}